\documentclass[review]{elsarticle}

\usepackage{lineno,hyperref}
\modulolinenumbers[5]

\usepackage{graphicx}
\usepackage[font = footnotesize, labelfont=bf]{subcaption}
\usepackage[font = footnotesize, labelfont=bf]{caption}

\usepackage[utf8]{inputenc}
\usepackage{natbib}
\usepackage{relsize}
\usepackage{xcolor}
\usepackage[margin=1.0in]{geometry}
\usepackage{fancyhdr}
\usepackage{amsmath}
\usepackage{bm}
\usepackage{setspace}
\usepackage{nicefrac}
\usepackage{indentfirst}
\usepackage{adjustbox}
\usepackage{enumitem}
\usepackage{algorithm}
\usepackage[noend]{algpseudocode}
\usepackage{subcaption}
\usepackage{cancel}

\usepackage{esint}
\numberwithin{equation}{section}

\begin{document}
	
	\begin{frontmatter}
		
		\title{\textbf{Implementation and performance analysis of efficient grid-free integral wall models in unstructured-grid LES solvers }}
		
		\author{Imran Hayat\fnref{myfootnote1}}
		\author{George Ilhwan Park\fnref{myfootnote2}}
		
		\address{Department of Mechanical Engineering and Applied Mechanics,
		
		University of Pennsylvania, 
		Philadelphia, Pennsylvania 19104, USA}
		\fntext[myfootnote1]{Ph.D. Candidate, Department of Mechanical Engineering and Applied Mechanics}
		\fntext[myfootnote2]{Assistant Professor, Department of Mechanical Engineering and Applied Mechanics}
		
		
		
		
		\begin{abstract}
			Two zonal wall-models based on integral form of the boundary layer differential equations, albeit with algebraic complexity, have been implemented in an unstructured-grid cell-centered finite-volume
            LES solver. The first model is a novel implementation of the ODE equilibrium wall model, where the velocity profile is expressed in the integral form using the constant shear-stress layer assumption and the integral is evaluated using a spectral quadrature method, resulting in a local and algebraic (grid-free) formulation. The second model, which closely follows the integral wall model of Yang \textit{et al.} (Phys. Fluids 27, 025112 (2015)), is based on the vertically-integrated thin-boundary-layer PDE along with a prescribed composite velocity profile in the wall-modeled region. Several numerical challenges unique to the implementation of these integral models in unstructured mesh environments, such as the exchange of wall quantities between wall faces and LES cells, and the computation of surface gradients, are identified and possible remedies are proposed. The performance of the wall models is assessed both in \textit{a priori} and  \textit{a posteriori} settings against the traditional finite-volume based ODE equilibrium wall model, showing a comparable computational cost for the integral wall model, and superior performance for the spectral implementation over the finite-volume based approach. Load imbalance among the processors in parallel simulations seems to severely degrade the parallel efficiency of finite-volume based ODE wall model, whereas the spectral implementation is remarkably agnostic to these effects.
		\end{abstract}
		
		
	\end{frontmatter}
	

	\section{Introduction }
	
	Different numerical approaches with varying degrees of accuracy and cost have been adopted over the past few decades to handle the multi-scale nature of turbulent flows. The most cost-effective and commonly employed technique in engineering applications, called Reynolds-averaged Navier-Stokes equations (RANS), tackles this problem by ensemble-averaging the Navier Stokes equations and thus modeling all the eddies present in the flow instead of resolving them. Clearly, the physics of turbulence, which may be critical to many problems (for instance, mixing flows and calculation of skin-friction coefficient in aerospace applications), is lost as a result of such simplistic modeling. On the other extreme, Direct Numerical Simulation (DNS) resolves all scales of motion present in the flow, providing the most accurate numerical solution. However, the cost of such simulations is prohibitive for any practical flows, especially at high Reynolds number. The number of grid points required for DNS is directly proportional to the ratio of largest to smallest scales of motion and since Reynolds number is a measure of this scale disparity, the required grid points scale with the Reynolds number according to $N^3 \propto Re_{L}^{9/4}$ \citep{piomelli02}.
    
    Large-Eddy simulation (LES) provides the middle ground between RANS and DNS, with lower cost than DNS and higher accuracy than RANS. In LES, the larger eddies typically containing more turbulent kinetic energy are resolved directly on the grid, whereas smaller eddies which tend to be universal and isotropic are modeled. LES solves filtered Navier-Stokes equations such that the transport of eddies larger than the size of filter is represented directly by the equations and the effect of smaller eddies is modeled through an artificial flux term called sub-grid scale (SGS) stress \citep{meneveau2000}, which is dynamically computed based on the LES solution, thus rendering LES a predictive method \citep{germano91,moin91}. The turbulent energy-representing characteristic of LES becomes its limitation near the wall in high Reynolds number wall-bounded flows. The grid resolution requirement of LES in the outer layer of turbulent boundary layer (TBL) is a linear function of Reynolds number \citep{choi12}, whereas in the inner layer, where smaller eddies account for most of the turbulent kinetic energy, the total number of grid points required scale as $O(Re_{\tau}^2)$ \citep{larsson16}. Therefore, to apply LES to a high Reynolds number TBL one would have to resolve the smaller scales near the wall, resulting in so-called Wall-Resolved LES (WRLES) which requires an enormous amount of grid points near the wall, thus eliminating the cost-effectiveness of LES over DNS \citep{bose18}. 
    
    This gives rise to the concept of Wall-Modeled LES (WMLES). In WMLES, the outer layer of TBL is resolved using a coarse LES, whereas the effect of inner layer dynamics is completely modeled. One common approach in the literature, called the wall-flux modeling, is to impose Neumann boundary at the wall by predicting the wall-fluxes using a wall-model, which takes inputs from the LES at some location in the inner layer and outputs the wall-flux to be imposed at the wall. Over the years, various models have been developed within the category of wall-flux modeling, ranging in complexity from algebraic to differential. The simplest wall-flux model called the algebraic equilibrium wall model relies on the law-of-the-wall to predict wall shear stress, thus making it devoid of most of the near-wall dynamics. In contrast, the two-layer zonal wall models solve either the full PDE for boundary layer or somewhat simplified ODE, typically with RANS parameterization. The ODE equilibrium wall-model \citep{wang02,bodart2012sensor,kawai12}, which is perhaps the most commonly employed wall-model for WMLES, ignores the non-equilibrium effects such as unsteadiness, non-linear advection and pressure gradient in favor of ease and low cost of implementation of the resulting ODE for diffusion. This ODE requires only wall-normal discretization, thus making this model local and essentially algebraic. However, this model still finds limited applicability in more practical and ubiquitous non-equilibrium flows. Alternatively, full unsteady three-dimensional RANS equations are solved employing a separate 3D grid for the wall-model \citep{park16jcp,park16prf}. Although this model incorporates most of the non-equilibrium effects near the wall, the associated cost for solving the wall model often approaches that of the main LES solution in the unstructured-grid solvers \citep{bose18}. Recently a new model called the Integral wall model, was formulated by Yang \textit{et al.} \citep{yang15}, which tries to strike a balance between the ODE equilibrium and the PDE wall models by considering the vertically integrated boundary layer equations while retaining all the non-equilibrium terms. The integration step is performed analytically  by assuming a general tractable form of the velocity profile in the inner layer based on the law-of-the-wall, thus reducing the PDE to an ODE. 
    
    Keeping in view the overarching theme of making WMLES more affordable in practical flow scenarios, which is one of the goals NASA's CFD Vision report \citep{slotnick2014cfd}, this paper details the implementation of efficient integral formulations of two zonal wall models namely, the ODE equilibrium wall model and the PDE non-equilibrium integral wall model in an unstructured-grid solver. A remark is in order regarding the scope of the \textit{integral wall model} terminology; thus far, the term \textit{integral wall model} has been used in the literature to refer to a model with vertically-integrated boundary layer PDE; however, in this paper we extend the scope of this term to include the spectral implementation of the ODE equilibrium wall model, since it involves the analytical integration of the constant-shear-stress equation. To make the distinction between ODE and PDE integral wall models clear, we refer to the original integral wall model by Yang \textit{et al.} \citep{yang15} as i\textit{ntegral nonequilibrium wall model} abbreviated as \textit{integral NEQWM}.
    
    Two points are worth noting here for the ODE equilibrium wall model. Firstly, despite ignoring the non-equilibrium effects, this model has still been shown to provide reasonable accuracy when deployed in non-equilibrium flows. This is attributed to the observation that even in strongly non-equilibrium flows the convective and pressure gradient terms approximately balance each other in the log-layer \citep{larsson16}. Secondly, despite its lower computational cost compared with more complex PDE models, the cost of solving a tridiagonal system (based on the finite-volumes (FV) solution to the ODE) is not negligible, typically accounting for 20-40 \% of the total LES cost in unstructured solvers \citep{bodart2012sensor,park17aiaa}. Furthermore, the cost overhead of the wall model is expected to aggravate in the parallel settings if an appropriate load-balancing strategy is not employed for the wall model solver while partitioning the domain. Conversely, alleviating the parallel efficiency degradation comes at the cost of additional implementation overhead. Therefore, it is worth investing more effort on improving the performance and ease of implementation of this commonly used wall model. To this end, we explore a spectral based solution of the ODE equilibrium model; specifically, the Gauss Quadrature (GQ) method is employed along with the integral form of the constant shear-stress equation in the inner layer to solve for the friction velocity (and consequently wall shear-stress). Essentially, this constitutes a grid-free approach, thus rendering the ODE equilibrium wall-model truly algebraic in complexity. We refer to this spectral implementation of the wall model as Gauss-Quadrature ODE Wall Model (GQ ODEWM) or simply GQWM, whereas the traditional finite-volumes based ODE wall model will be abbreviated as FVWM from this point onward.
    
    On the other hand, the integral NEQWM with its already proven cost-effectiveness over differential-complexity models, is mainly approached in this study from the perspective of the strategy and the challenges associated with its implementation in an unstructured-grid cell-centered finite-volume LES solver. It should be noted that the original formulation of this model by Yang \textit{et al.} \citep{yang15} was implemented in a structured-grid finite-difference/spectral solver, and a later extension of the model to the compressible flows by Catchirayer \textit{et al.} \citep{catchirayer18} employed a finite-volume structured-grid solver. This makes the present study the first effort in the extension of this method to more versatile unstructured-grid, finite-volume LES solvers, which are frequently and preferably used in the industry and the research community alike due to their ability to handle complex geometries. Indeed, we observe unique challenges in the present study which were not encountered in the structured-grid environment. Furthermore, both the previous versions of integral NEQWM were limited to essentially one-dimensional wall shear stress (for the former, this resulted from an inconsistent choice of the viscous sublayer profile, and for the latter, this was done by design). 
    
    In this study, we discuss the implementation overhead of integral NEQWM in the unstructured-grid environment in detail and propose possible strategies to overcome these challenges. We present a modification to the original integral NEQWM formulation in \citep{yang15} to rectify the inconsistent near-wall asymptotic behavior of the velocity. A brief validation study is conducted for the two integral wall models using the channel flow. Performance characteristics of the spectral-based ODE wall model are discussed in detail both in \textit{a priori} and  \textit{a posteriori} settings. A note-worthy caveat in this study is that currently only the incompressible formulation of both these wall models have been implemented. This limitation will be addressed in future studies.
	
	The paper is organized as follows: \S \ref{sec:govern_eqn} details the governing equations used in the LES and the wall models. In \S \ref{sec:discretize_n_soln}, discretization of the wall model equations and their solution method is provided. In \S \ref{sec:WM_implementation}, implementation details of coupling between wall-model and LES solver, the communication algorithms as well as various numerical challenges in implementing the wall models are discussed. Performance characteristics of the wall models in \textit{a priori} and \textit{a posteriori} settings are analyzed in \S \ref{sec:performance_validation}, while \S \ref{sec:conclusion} provides concluding remarks.
	
	\section{Governing Equations }\label{sec:govern_eqn}

    The generic governing equations for the LES solver and the wall-model are the compressible filtered or ensembled-averaged Navier Stokes equations as shown in the conservative form below,
  
 	\begin{equation}
	\label{eqn:NS_continuity}
	\frac{\partial {\rho}}{\partial t}+\frac{\partial \rho {u}_{j} }{\partial x_{j}}=0
	\end{equation}   
    
	\begin{equation}
	\label{eqn:NS_momentum}
	\frac{\partial \rho{u}_{i}}{\partial t}+\frac{\partial \rho u_{i} {u}_{j} }{\partial x_{j}} + \frac{\partial p}{\partial x_{i}}=\frac{\partial \tau_{ij}}{\partial x_{j}}
	\end{equation}    

	\begin{equation}
	\label{eqn:NS_energy}
	\frac{\partial \rho E}{\partial t}+\frac{\partial (\rho E + p) {u}_{j} }{\partial x_{j}} =\frac{\partial \tau_{ij} u_{i}}{\partial x_{j}} - \frac{\partial q_{j}}{\partial x_{j}}
	\end{equation}   

    where, $\rho$, $u_{i}$, $p$, and $E$ represent the density, velocity, pressure and total energy, respectively. 
    Next we present the specific form of the governing equations for each of the wall models considered. It should be reiterated here that the equations are presented for the incompressible formulation. The thin boundary-layer momentum equation for a general wall-flux model is given as:
	
	\begin{equation}
	\label{eqn:WM}
	\frac{\partial {u}_{i}}{\partial t}+\frac{\partial u_{i} {u}_{j} }{\partial x_{j}}+\frac{1}{\rho} \frac{\partial p}{\partial x_{i}}=\frac{\partial}{\partial y}\left[\left(\nu+\nu_{t}\right) \frac{\partial {u}_{i}}{\partial y}\right],  \; \; \; \;  i=1,3,  \; \; \; \; j=1,2,3
	\end{equation}
	
	\noindent where $u_{i}$ is the local wall-parallel velocity in streamwise or spanwise ($x_{1}$ or $x_{3}$) direction, $y$ is the local wall-normal coordinate, $\nu$ is the kinematic viscosity and $\nu_{t}$ is the eddy-viscosity. All the flow variables in the above equation are ensemble-averaged quantities.
    
    \subsection{ODE equilibrium wall model}

	In the ODE equilibrium wall-model, the non-equilibrium effects in the momentum equation \ref{eqn:WM} are neglected, rendering left-hand-side equal to zero and reducing the PDE to an ODE,
	
	\begin{equation}
	\label{eqn:ODEWM}
	\frac{\partial}{\partial y}\left[\left(\nu+\nu_{t}\right) \frac{\partial u}{\partial y}\right]=0,
	\end{equation}
	
	\noindent where $u$ is the tangential velocity parallel to the wall. The terms in the square brackets constitute the total shear stress $\tau$. Note that when integrated, equation \ref{eqn:ODEWM} leads to a constant shear-stress across the wall-model region. Equation \ref{eqn:ODEWM} is typically solved using the finite-volume method on a separate wall-normal grid which is embedded on the LES grid between the wall and the matching-location, along with the no-slip boundary condition at the wall and the wall-parallel LES velocity at the matching location, to obtain wall-parallel velocity $u$ in each cell of the wall model grid. The wall shear stress can then be obtained from the resulting velocity profile using the definition of shear stress at the wall.

	\begin{equation} 
	\label{eqn:wallstress1}
	\tau_{w}=\mu \left.\frac{du}{dy}\right\vert_{w}  \approx \mu \frac{u_{1}-0}{\Delta y_{1}}
	\end{equation}
	
	\noindent where the velocity gradient at the wall is approximated using the first off-wall cell value $u_{1}$ and half cell-height $\Delta y_{1}$ in the wall-model grid along with the no-slip condition.

	Note that the above method requires an iterative solution of a tri-diagonal system through matrix inversion at each wall face until the value of $\tau_{w}$ converges. The size of the tri-diagonal system is equal to the number of points in the ODE grid in the wall-normal direction. This makes the ODE wall-model implementation somewhat computationally intensive, especially for a large number of wall faces. 
	
	
	\subsubsection{Gauss quadrature based ODE equilibrium wall-model}
	An alternative form of the ODE Equilibrium model in the integral form is presented below, which retains the complete inner-layer profile down to the wall. This reformulation makes the model amenable to Gauss-quadrature based (grid-free) methods, which could prove to be more efficient than the finite-volume approach, especially when the wall-model requires a larger grid size. The reformulation is possible by making use of the Van Driest damping function, as detailed in \citep{pope} and summarized below for completeness.
	
	
	The total shear stress in the inner layer can be expressed as,
	
	\begin{equation}
	\tau(y) / \rho = \nu \frac{\partial u}{\partial y}+\nu_{\mathrm{t}} \frac{\partial u}{\partial y} \\
	\end{equation}
	
	\begin{equation}\label{eqn:totalstress}
	\tau(y) / \rho = \nu \frac{\partial u}{\partial y}+\ell_{\mathrm{m}}^{2}\left(\frac{\partial u}{\partial y}\right)^{2}   
	\end{equation}
	
	\noindent where $l_{m}$ is the mixing length. Equation \ref{eqn:totalstress} assumes $\frac{\partial u}{\partial y}>0$ in the inner layer, which limits this model to attached flows. In this equation, the eddy-viscosity $\nu_{t}$ has been approximated by the mixing-length hypothesis,
	
	\begin{equation}\label{eqn:Prandtl}
	\nu_{\mathrm{t}}=\ell_{\mathrm{m}}^{2}\left|\frac{\mathrm{d} u}{\mathrm{d} y}\right|.
	\end{equation}
	
	\noindent Normalizing equation \ref{eqn:totalstress} by viscous scales and restating the equation in wall-units ($y^{+}=y/\delta_{\nu}$, $u^{+}=u/u_{\tau}$, $l_{m}^{+}=l_{m}/\delta_{\nu}$) we get,
	
	\begin{equation}\label{eqn:totalstress_norm}	
	\frac{\tau}{\tau_{\mathrm{w}}}=\frac{\partial u^{+}}{\partial y^{+}}+\left(\ell_{\mathrm{m}}^{+} \frac{\partial u^{+}}{\partial y^{+}}\right)^{2}.
	\end{equation}
	
	\noindent We use the definition of mixing length $l_{m}$ given by van Driest \citep{van56}, which closely approximates the profiles in the log-layer, the buffer-layer and the viscous sublayer for an equilibrium boundary layer,
	
	\begin{equation}\label{eqn:lm_VD}
	\ell_{\mathrm{m}}^{+}=\kappa y^{+} \left[1-\mathrm{exp}\left(-y^{+}/A^{+}\right)\right],
	\end{equation}
	
	\noindent where $A^{+}$ is a constant with value $A^{+}=26$.
	
	The constancy of total shear stress as given by equation \ref{eqn:ODEWM} results in $\tau=\tau_{\mathrm{w}}$ along the entire height of the wall-model region. Therefore, the left-hand-side of equation \ref{eqn:totalstress_norm} is equal to unity. This is a quadratic equation in $\frac{\partial u^{+}}{\partial y^{+}}$, which can be solved as,
	
	\begin{equation}\label{eqn:dudy}
	\frac{\partial u^{+}}{\partial y^{+}}=\frac{2}{1+\left[1+\left(4\right)\left(\ell_{\mathrm{m}}^{+}\right)^{2}\right]^{1 / 2}}.
	\end{equation}

	
	\noindent Integrating the dimensional form of equation \ref{eqn:dudy} from the wall up to an arbitrary height $y$ using no-slip condition at the wall gives the integral expression for velocity $u(y)$.
	
	\begin{equation}\label{eqn:u_integral}
	u(y)=\int_{0}^{y} \frac{2 u_{\tau}^{2}}{\nu} \frac{1}{1+\left[1+4\left(\ell_{m}^{+}(y^{\prime})\right)^{2}\right]^{\frac{1}{2}}} d y^{\prime}.
	\end{equation}
	
	\noindent Equation \ref{eqn:u_integral} provides the complete continuous velocity profile in the inner layer, unlike the composite profiles which consider piece-wise variation in log-layer and viscous sublayer. Note that by imposing the velocity at the matching location $u(y=h_{wm})=U_{LES}$ in equation \ref{eqn:u_integral}, we get an integral equation with $u_{\tau}$ as the only unknown,
	
	\begin{equation}\label{eqn:utau_equation}
	U_{LES} - \int_{0}^{h_{wm}} \frac{2 u_{\tau}^{2}}{\nu} \frac{1}{1+\left[1+4\left(\ell_{m}^{+}(y^{\prime})\right)^{2}\right]^{\frac{1}{2}}} d y^{\prime} = 0 .
	\end{equation}	
 
	\noindent Equation \ref{eqn:utau_equation}, which is an integral statement of the ODE equilibrium wall model, allows us to circumvent the need for a grid to solve this model, by a employing grid-free quadrature method to evaluate the integral.

	\subsection{Integral NEQWM}
	
    The integral NEQWM includes the pressure-gradient and acceleration terms in equation \ref{eqn:WM}. The governing equations are RANS-type equations with thin boundary-layer approximation applied to them. Most PDE wall-models numerically integrate differential equations of the form \ref{eqn:WM} along $y$-direction using a vertically refined mesh. The integral wall-model circumvents this step by assuming a parameterized form of the velocity profile in the $y$-direction, based on the known scaling laws in the inner-layer, and performs a vertical integration of the equations analytically. The resulting equation is marched in time synchronously with the LES. This approach is similar to Von-Karman-Pohlhausen integral method \citep{karman21, pohlhausen21}, which assumes a polynomial form of velocity profile to analytically integrate the momentum equations, thus relating the wall-stress with the prescribed velocity value at some distance from the wall. This method is computationally less intensive compared to other non-equilibrium PDE wall-models because no grid is required in the vertical direction; furthermore, the wall-parallel spatial gradients prescribed at the matching location are computed from the known neighborhood values at the previous time step. In this sense the integral wall-model is local and algebraic. Therefore, the model retains all the relevant near-wall physics while maintaining the complexity of an equilibrium wall-model. The governing equations for this model are given below.
    
    \begin{equation}\label{eqn:iWMCont} 
    \frac{\partial u}{\partial x}+\frac{\partial v}{\partial y}+\frac{\partial w}{\partial z}=0
    \end{equation}
    
    \begin{equation}\label{eqn:iWMRANSy}
    \frac{\partial p}{\partial y}=0
    \end{equation}
    
    \begin{equation}\label{eqn:iWMRANSx}
    \frac{\partial u}{\partial t} +                                   \frac{\partial u \,u}{\partial x} +                     \frac{\partial u \, v}{\partial y} +                      \frac{\partial u \, w}{\partial z} = 
    -\frac{1}{\rho} \frac{\partial p}{\partial x} 
    +\frac{\partial}{\partial y}\left[\left(v+v_{T}\right) \frac{\partial u}{\partial y}\right]
    \end{equation}

    \begin{equation}\label{eqn:iWMRANSz}
    \frac{\partial w}{\partial t} + \frac{\partial w\, u}{\partial x} + \frac{\partial w\, v}{\partial y} + \frac{\partial w\, w}{\partial z} = 
    -\frac{1}{\rho} \frac{\partial p}{\partial z} 
    +\frac{\partial}{\partial y}\left[\left(v+v_{T}\right)   \frac{\partial w}{\partial y}\right]
    \end{equation}

    
    \noindent Equation \ref{eqn:iWMCont} is the mean continuity and equations \ref{eqn:iWMRANSy} to \ref{eqn:iWMRANSz} are the expanded form of equation \ref{eqn:WM}.  Equations \ref{eqn:iWMCont} to \ref{eqn:iWMRANSz} are integrated vertically from $y=0$ to the matching location $y=h_{wm}$ to obtain integral momentum equations in $x$ and $z$ directions, which are eventually expressed in terms of the parameters of the assumed velocity profile. These two equations are solved along with a set of other physical and matching conditions to obtain the unknown parameters of the profile. This process is described below.
    
    For the purpose of delineating the steps of this method in a tractable way, here we show the simplified formulation for a two-dimensional flow in local $x-y$ (streamwise-wallnormal) plane. However, it should be noted that the wall model in this study has been implemented for full 3D boundary layer equations as expressed in equations \ref{eqn:iWMCont} to \ref{eqn:iWMRANSz}. For complete details of the 3D formulation, the reader is referred to \citep{yang15}, which was closely followed (in terms of formulation and notations) in the present implementation of the integral NEQWM. With the 2D assumption, integrating equation \ref{eqn:iWMCont} vertically from $y=0$ to $y=h_{wm}$ and noting that $\left. v\right|_{y=0}=0$, we get,

	\begin{equation}\label{eqn:iWMCont2}
	\left.  v\right|_{y=h_{wm}}=-\frac{\partial}{\partial x} \int_{0}^{h_{wm}}  u d y.
	\end{equation}

    \noindent Integrating equation \ref{eqn:iWMRANSx} vertically and noting that from equation \ref{eqn:iWMRANSy} that $ p= p(x)$ ,
	
	\begin{equation}\label{eqn:iWMmomxint}
	\frac{\partial}{\partial t}\int_{0}^{h_{wm}} u dy + 
	\frac{\partial }{\partial x} \int_{0}^{h_{wm}} u^2 dy +
	\left. u \,  v\right|_{y=h_{wm}}-\left. u \,  v\right|_{y=0} = 
	\frac{1}{\rho} \left[-\frac{\partial p}{\partial x} h_{wm} 
	+\tau_{h_{wm}} - \tau_{w}\right].
	\end{equation}

    \noindent Substituting $\left. v\right|_{y=h_{wm}}$ from equation \ref{eqn:iWMCont2} in equation \ref{eqn:iWMmomxint} and noting that $\left. u\right|_{y=0}=0$ and $\left. u\right|_{y=h_{wm}}=U_{LES}$,
	
	\begin{equation}
	\frac{\partial}{\partial t}\int_{0}^{h_{wm}} u dy + 
	\frac{\partial }{\partial x} \int_{0}^{h_{wm}} u^2 dy -
	U_{LES} \,\frac{\partial}{\partial x} \int_{0}^{h_{wm}}  u d y = 
	\frac{1}{\rho} \left[-\frac{\partial p}{\partial x} h_{wm} 
	+\tau_{h_{wm}} - \tau_{w}\right],
	\end{equation}
	

    \noindent or,
	
	\begin{equation}\label{eqn:iWMx}
	\frac{\partial L_{x}} {\partial t} + 
	\frac{\partial L_{xx}}{\partial x} -
	U_{LES} \,\frac{\partial L_{x}}{\partial x}= 
	\frac{1}{\rho} \left[-\frac{\partial p}{\partial x} h_{wm} 
	+\tau_{h_{wm}} - \tau_{w}\right],
	\end{equation}
	
    \noindent where $	L_{x} \equiv \int_{0}^{h_{wm}} u dy$ , 
	$L_{xx}\equiv \int_{0}^{h_{wm}} u^2 dy$ ,	$\tau_{h_{wm}} = \left(\mu+\mu_{T}\right) \left.\frac{\partial u}{\partial y}\right|_{y=h_{wm}} $ and $ \tau_{w} = \mu \left.\frac{\partial u}{\partial y}\right|_{y=0}$. Equation \ref{eqn:iWMx} is the integral form of boundary layer $x$-momentum equation. Note that $U_{LES}$ in this equation is the time-filtered velocity from LES at the matching location $y=h_{wm}$. The terms $L_{x}$ and $L_{xx}$ (hereinafter called the \textit{integral terms of the wall-model} or more concisely the \textit{integral terms}) are unclosed. In order to close the integral terms, a functional form of the velocity profile is needed. Based on the scaling laws in the turbulent boundary layer, Yang \textit{et al.} \citep{yang15} proposed the following parametric form of the profile for the the linear viscous-sublayer and the logarithmic-layer (log-layer),
    
    \begin{equation}\label{eqn:iWMprofile}
    \begin{array}{ll}
     u = u_{\tau} \frac{y}{\delta_{\nu}} =  \frac{u_{\tau}^2}{\nu}y, & 0 \leq y \leq \delta_{i} \\
     u = u_{\tau}\left[\frac{1}{\kappa} \log \frac{y}{h_{wm}}+C\right]+u_{\tau} A \frac{y}{h_{wm}}, & \delta_{i}<y \leq h_{wm},
    \end{array}
    \end{equation}
    
    \noindent where the four unknown parameters are defined as: $\delta_{i}$ is the height of viscous sublayer, $u_{\tau}= \sqrt{\tau_{w}/\rho}$ is the friction velocity, coefficient $C$ ensures C$^{0}$ continuity of the profile at $\delta_{i}$, and coefficient $A$ determines the linear variation from the log-law as a result of non-equilibrium effects (pressure gradient and advection); $\delta_{\nu}=\nu/u_{\tau}$ is the viscous length-scale. Note that for the 3D formulation consisting of two wall-parallel velocity components, a total of 8 unknown parameters ($u_{\tau},u_{\tau,x},u_{\tau,z},A_{x},A_{z},C_{x},C_{z},\delta_{i}$) would be needed to describe the composite profiles in $x$ and $z$.
    
    \noindent The Integral terms can now be evaluated as follows:
    
    \begin{equation}\label{eqn:Lx1}
    \begin{aligned}
    L_{x} &= \int_{0}^{\delta_{i}}  u dy + \int_{\delta_{i}}^{h_{wm}}  u dy \\
    L_{x} &=\frac{1}{2} u_{\tau} \frac{\delta_{i}^{2}}{\delta_{\nu}}+u_{\tau} h_{wm}\left[\frac{1}{2} A\left(1-\frac{\delta_{i}^{2}}{h_{wm}^{2}}\right)-\frac{1}{\kappa}+\left(1-\frac{\delta_{i}}{h_{wm}}\right)\left(C+\frac{1}{\kappa} \log \frac{\delta_{i}}{h_{wm}}\right)\right]
    \end{aligned}
    \end{equation}
    
    \begin{equation}\label{eqn:Lxx1}
    \begin{aligned}
    L_{xx} = \int_{0}^{\delta_{i}}  u^2 dy &+ \int_{\delta_{i}}^{h_{wm}}  u^2 dy \\
    L_{xx}=\frac{1}{3} u_{\tau}^{2} \frac{\delta_{i}^{3}}{\delta_{v}^{2}}+&u_{\tau}^{2} h_{wm}\left[\left(C-\frac{1}{\kappa}\right)^{2}-\frac{\delta_{i}}{h_{wm}}\left(C-\frac{1}{\kappa}+\frac{1}{\kappa} \log \frac{\delta_{i}}{h_{wm}}\right)^{2}+\frac{1}{\kappa^{2}}\left(1-\frac{\delta_{i}}{h_{wm}}\right)\right. \\
    &+A\left(C-\frac{1}{2 \kappa}\right)\left(1-\frac{\delta_{i}^{2}}{h_{wm}^{2}}\right)-\frac{A}{\kappa} \frac{\delta_{i}^{2}}{h_{wm}^{2}} \log \frac{\delta_{i}}{h_{wm}}+\frac{1}{3} A^{2}\left(1-\frac{\delta_{i}^{3}}{h_{wm}^{3}}\right).
    \end{aligned}
    \end{equation}
    
    \noindent Furthermore, the shear stresses in equation \ref{eqn:iWMx} can be expressed as:
    
    \begin{equation}\label{eqn:tau1}
    \begin{aligned}
    \tau_{h_{wm}} &= \rho \left(\mu+\left.\mu_{t}\right|_{y=h_{wm}}\right) \frac{u_{\tau}}{h_{wm}}\left(\frac{1}{\kappa}+A\right) \\
    \tau_{w} &= \rho u_{\tau}^2
    \end{aligned}
    \end{equation}
    
    \noindent where $\mu_{t}$ is given by equation \ref{eqn:Prandtl}	

    \subsection{Modifications to the original integral NEQWM}

    In the original 3D integral NEQWM formulation in \citep{yang15}, the assumed velocity profiles for the two wall-parallel velocity components in the viscous sublayer are given by equation (C5) in \citep{yang15} as,
    
    \begin{eqnarray}\label{eqn:sublayer_inconsistent}
    u = u_{\tau,x} \frac{y}{\delta_{\nu}} \nonumber \\
    w = u_{\tau,z} \frac{y}{\delta_{\nu}}.
    \end{eqnarray}       
    
    \noindent With the above formulation, we observed that the wall-stress predictions of the wall model were highly sensitive to the choice of the local $x$-$z$ coordinates. This anomaly was traced to the inconsistent asymptotic behavior of the sublayer velocity profile near the wall. The desired asymptotic behavior of the velocity near the wall is given by the Taylor series expansion as, 
    
    \begin{equation}\label{eqn:true_asypmtotic}
    u = \cancelto{\textrm{no-slip}}{u(0)} \qquad \quad+ \left.\frac{\partial u}{\partial y}\right|_{w} y + \mathcal{O}(y^2)  = \left(\frac{\tau_{w,x}}{\mu}\right) y,
    \end{equation}     
    
    \noindent whereas with the formulation \ref{eqn:sublayer_inconsistent}, the following asymptotic behavior is obtained (see \ref{appendix:A} for derivation),
    
    \begin{equation}\label{eqn:false_asypmtotic}
    u = \left(\frac{\tau_{w,x}}{\mu}\right) y \frac{1}{\sqrt{\cos^2\theta + \cos \theta \sin \theta}},
    \end{equation}

    \noindent where $\theta = \arctan \left(\frac{\tau_{w,z}}{\tau_{w,x}}\right) $. The extra factor at the end of equation \ref{eqn:false_asypmtotic} renders the original sublayer formulation inconsistent. In our current formulation, we have modified the assumed viscous sublayer profile as follows, 
    
    \begin{eqnarray}\label{eqn:sublayer_modified}
    u = sign(u_{\tau,x}) \frac{u_{\tau,x}^{2}}{u_{\tau}} \frac{y}{\delta_{\nu}} \nonumber \\
    w = sign(u_{\tau,z}) \frac{u_{\tau,z}^{2}}{u_{\tau}} \frac{y}{\delta_{\nu}}.
    \end{eqnarray}
    
    \noindent The above choice of assumed profile ensures a consistent near-wall asymptotic behavior given by equation \ref{eqn:true_asypmtotic}, as shown in \ref{appendix:A}. 
	
	\section{Discretization and solution method}\label{sec:discretize_n_soln}
	
	\subsection{Solver information}
	
	Both the wall models considered in this study have been integrated as sub-routines into CharLES, a cell-centered unstructured finite-volume compressible LES solver developed by Cascade Technologies Inc. The solver uses a second-order central scheme for spatial discretization and an explicit third-order Runge-Kutta (RK3) scheme for time advancement (opposed to the explicit Euler scheme used in the integral wall model). The code extensively employs the object-oriented programming (OOP) structure of C++ and also supports parallelism through Message Passing Interface (MPI). Many features of this code, including the main LES solver, have been optimized for load-balancing across multiple processors. However, as shown in previous studies, achieving parallelism across multiple routines is not always trivial, especially when the modules communicate with each other at every time-step, which is the case for wall models \citep{park16jcp}. This further necessitates the development of grid-free approaches discussed in this study in order to circumvent the need for grid-partitioning and parallelism within wall model routines. Details regarding other aspects of the flow solver can be found in \citep{park16jcp, khalighi11}

	\subsection{Discretization of Gauss-quadrature based ODE equilibrium wall model}
	
	The integral in equation \ref{eqn:utau_equation} is evaluated using the Gauss-Lobatto-Legendre (GLL) quadrature method. A simple Gauss-Legendre quadrature excludes the boundary points, therefore the GLL is chosen to include the boundary conditions in the evaluation of the integral. The $n$-point GLL quadrature formula is given by,
	
	\begin{equation}\label{eqn:GLL}
	\int_{-1}^{1} f(\xi) d \xi=\sum_{i=0}^{n-1} w_{i} f\left(\xi_{i}\right)+R,
	\end{equation}
	
	\noindent where $f(\xi)$ is the integrand, $n$ is total number of quadrature points, $\xi_{i}$ are the quadrature points or abscissae and $w_{i}$ are the weights of the GLL quadrature, $R$ is the error in approximation of the integral, with $R=0$ giving the exact integral. For optimal accuracy using GLL, the quadrature points are taken as the boundary points $\xi=\pm1$ and the zeros of the Jacobi polynomial $P_{n-2}^{1,1}(\xi)$ (or equivalently the zeros of $\frac{d}{d \xi}\left(L_{n-1}(\xi)\right)$, where $L_{n-1}(\xi)$ is the Legendre polynomial); and the weights are given by $	w_{i}=\frac{2}{n(n-1)\left[L_{n-1}\left(\xi_{i}\right)\right]^{2}}$. The formula for GLL quadrature gives the exact integral when polynomial order of the integrand is $N \leq 2n-3$ i.e. $	R=0 \,\,\,\,\, \forall \,\,\,\,\, f(\xi) \in \mathcal{P}_{2n-3}([-1,1])$, where $\mathcal{P}_{2n-3}([-1,1])$ denotes the space of polynomials up to order $(2n-3)$ with the domain $[-1,1]$.
	
	To evaluate the integral in equation \ref{eqn:utau_equation} using the quadrature formula given by equation \ref{eqn:GLL}, a change of variable $\xi \to y$ from the Legendre domain $[-1,1]$ to the physical domain $[0, h_{wm}]$ must be applied. This can be done using the following linear transformation, which maps the end points of the Legendre interval ($\xi = -1$ and $\xi = 1$) to the physical boundary points of the wall model ($y = 0$ and $y = h_{wm}$, respectively). 
	
	\begin{equation}\label{eqn:lin_transf}
	y=\frac{h_{wm}}{2} (1 + \xi).
	\end{equation}
	
	\noindent The corresponding quadrature formula in the physical domain for the integral in equation \ref{eqn:utau_equation} is given by,
	
	\begin{equation}\label{eqn:GLL_u}
	\int_{0}^{h_{wm}} I(y) d y = \frac{h_{wm}}{2} \sum_{i=0}^{n-1} w_{i} I\left(\frac{h_{wm}}{2} + \frac{h_{wm}}{2}\xi_{i}\right),
	\end{equation}	
	
	\noindent where, \qquad \qquad \qquad \qquad \qquad \qquad $I(y)=\frac{2 u_{\tau}^{2}}{\nu} \frac{1}{1+\left[1+4\left(\ell_{m}^{+}(y)\right)^{2}\right]^{\frac{1}{2}}}.$
	
    \bigskip    
    
	Equation \ref{eqn:GLL_u} can be solved iteratively for $u_{\tau}$ by using a shooting method (e.g. secant method) to evaluate the integral. Note that unlike the traditional ODE model, which solves for the entire $u(y)$ profile in the wall-model region, the proposed method only solves for $u_{\tau}$. However, the entire velocity profile can be constructed from equation \ref{eqn:u_integral} using the predicted value of $u_{\tau}$.	
	
	Note that the linear transformation in equation \ref{eqn:lin_transf} maps the nodes to the physical domain such that most of the nodes are clustered close to the two end points of the domain ($y = 0$ and $y = h_{wm}$). However, a closer inspection shows that the integrand term varies more drastically closer to the wall than away from it. Therefore, to make the integral evaluation more efficient (i.e. using lesser number of quadrature points to approximate the integrand), a non-linear transformation of the following form can be employed, which maps more nodes closer to the wall and fewer nodes away from it.
	
	\begin{equation}\label{eqn:nonlin_transf}
	y=h_{wm} \frac{e^{\xi + 1}-1}{e^{2}-1}.
	\end{equation}	
	
	\noindent The corresponding quadrature formula is given as,
	
	\begin{equation}\label{eqn:GLL_u_nonlin}
	\int_{0}^{h_{wm}} I(y) d y = \frac{h_{wm}}{e^{2}-1} \sum_{i=0}^{n-1} w_{i} e^{\xi_{i} + 1} I\left(h_{wm} \frac{e^{\xi_{i} + 1}-1}{e^{2}-1}\right).
	\end{equation}

    \noindent The cost advantage of using a non-linear transformation will become evident in \S \ref{sec:GQ_perform}.

	\subsection{Discretization and solution method of integral NEQWM}\label{sec:IWM_discretized}
	
	Below, we delineate the steps to obtain the four unknown parameters in the velocity profile in equation (\ref{eqn:iWMprofile}).

    1. Impose LES velocity at the matching location:
    
    \begin{equation}\label{eqn:cond1}
    u_\tau\left[C+A\right] = U_{LES}
    \end{equation}
    
    2. Impose continuity of the profile at the interface of viscous sublayer and log-layer $y=\delta_{i}$:
    
    \begin{equation}\label{eqn:cond2}
    \frac{u_{\tau}^2}{\nu}\delta_{i}= u_{\tau}\left[\frac{1}{\kappa} \log \frac{\delta_{i}}{h_{wm}}+C\right]+u_{\tau} A \frac{\delta_{i}}{h_{wm}}
    \end{equation}
    
    3. Determine $\delta_{i}$ from the continuity of viscous-sublayer and standard log-law at $y=\delta_{i}$ with constants $\kappa=0.4$ and $B=5$. 
    
    \begin{equation}\label{eqn:cond3}
    \begin{aligned}
    \delta_{i}^{+}&= \left[\frac{1}{\kappa} \log \delta_{i}^{+}+B\right] \\ \Rightarrow \, \delta_{i}^{+}&=11
    \end{aligned}
    \end{equation}
    
    \noindent where $\delta_{i}^{+}=\delta_{i}/\delta_{\nu}$ is the non-dimensional viscous sublayer height in wall-units.
    
    4. The four equations \ref{eqn:iWMx} and \ref{eqn:cond1}$-$\ref{eqn:cond3} can now be solved simultaneously to obtain the 4 unknown parameters. Note that in equation \ref{eqn:iWMx} the integral terms and the shear stresses are completely expressible in terms of the unknown parameters of the profile, as evident from equations \ref{eqn:Lx1} to \ref{eqn:tau1}.
    
    
    5. Equation \ref{eqn:iWMx} is a PDE in $x$ and $t$, that can be solved numerically. The time integration is performed using explicit Euler scheme such that all the spatial derivative terms and shear-stresses are evaluated at time-step $n$ and therefore are known quantities for determining $L_{x}$ at time step $n+1$. The time-discretized form of equation \ref{eqn:iWMx} is as follows,
    
    \begin{equation}\label{eqn:iWMx_discretized}
    L_{x}^{(n+1)} = L_{x}^{(n)} + \Delta t \left[-\frac{\partial L_{xx}}{\partial x} + U_{LES} \,\frac{\partial L_{x}}{\partial x} 
    +\frac{1}{\rho} \left(-\frac{\partial p}{\partial x} h_{wm} 
    +\tau_{h_{wm}} - \tau_{w}\right)\right]^{(n)}.
    \end{equation}
    
     Equation \ref{eqn:iWMx_discretized} being non-linear, is solved iteratively using Newton-Raphson method up to a user-specified error tolerance. In solving this equation, $U_{LES}$ and pressure gradient are the only inputs from the resolved LES solution at the specified matching location $h_{wm}$. The spatial derivatives in this equation can be approximated using a bi-linear interpolation method, which for a structured grid reduces to the second-order central differencing if the neighbouring wall-face values are used. Since $L_{x}^{(n)}$ and $L_{xx}^{(n)}$ are known for all wall-faces, the derivative terms $\frac{\partial}{\partial x}(.)$ could in principle be readily approximated using a difference formula of the form,
    
    \begin{equation}\label{eqn:FD}
    \frac{\partial \Phi_{i}^{(n)}}{\partial x} = \frac{\Phi_{i+1}^{(n)}-\Phi_{i-1}^{(n)}}{2 \Delta x},
    \end{equation}

    \noindent where $\Phi$ is one of the scalars: $L_{x}$ , $L_{xx}$ or $ p$. However, this approach for computing spatial gradients proposed in the original study by Yang \textit{et al.} \citep{yang15} is not well-suited to the current study for two reasons. Firstly, CharLES solver being an unstructured-grid solver, does not have geometrically ordered connectivity between cells to employ simplistic reductions like equation \ref{eqn:FD}; instead it computes gradients using the Green-Gauss second-order reconstruction, that uses values from all immediate CV neighbors of a cell \citep{park16jcp}. Secondly the full 3D integral equations in $x$ and $z$ (see equation \ref{eqn:iWMx_3D}) require coordinate transformations of the wall quantities (integral terms) from local to global coordinate system and vice-versa, which imposes an additional implementation overhead. A detailed discussion on these numerical aspects is provided in \S \ref{sec:surface_grad} and \S \ref{sec:surface_grad_prob}.

	\section{Numerical aspects of wall-model implementation in unstructured LES solver}\label{sec:WM_implementation}
	
	To facilitate interested readers to readily test, modify, and adapt to their own solvers the wall models discussed in this study, we provide simple modular MATLAB implementations of all three wall models discussed in this study (GQ-EQWM, FV-EQWM, and integral NEQWM) in our GitHub repository \url{https://www.github.com/imranhayat29/Wall-Models-for-LES}. However, it must be cautioned that these routines were developed mainly with \textit{apriori} validation in mind and as such, some of the implementation aspects (such as those related to the exchange of data between LES mesh and wall model, and the computation of surface gradients in unstructured solvers) may have been oversimplified/avoided in these routines. These implementation challenges become truly evident only when the wall models are deployed in actual solvers and coupled with the LES mesh.
	
	\subsection{Implementation of coupled WMLES solver}\label{sec:coupled_WMLES}
 
\begin{figure}[t]
			\includegraphics[trim=0 0 250 0,clip,width=\linewidth]{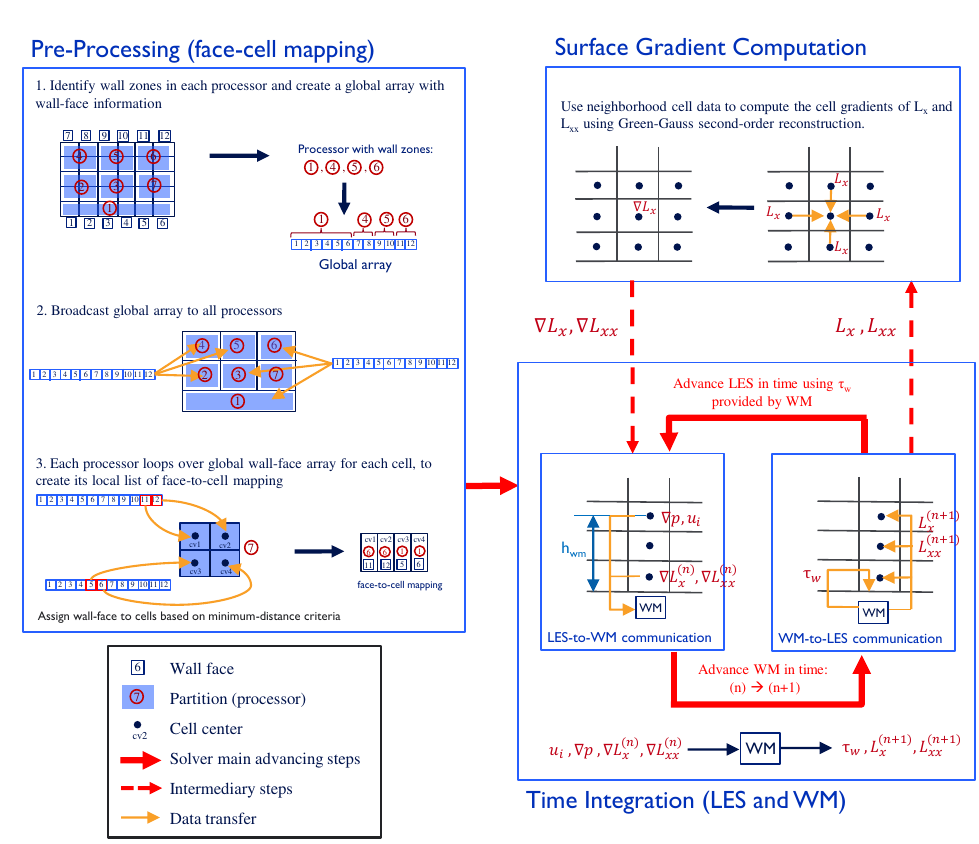} 
		
        \caption[]
        {\small Schematic of integral NEQWM routine coupled with the LES solver.} 
        
        \label{fig:algorithm_schematic}
	\end{figure}

	A generic difference between the two models considered in this study is noted here from the implementation standpoint. Since the ODE Equilibrium wall model is a one-dimensional and truly local (i.e. it does not require information from the neighboring wall-faces), for a given wall-face, its implementation and coupling to the LES solver is restricted to a single point in the LES grid (the matching location) and does not require communication protocols with the rest of the LES grid. On the contrary, integral NEQWM requires data from the neighbouring wall-faces to compute surface gradients, which necessitates the development of a proper algorithm with well-defined communication protocols for this method, especially given the parallelism and partitioning in the LES solver. With this in mind, a general algorithm is provided for the integral NEQWM in algorithm \ref{algo:WM}. It is however noted that most of the generic steps in this algorithm are still applicable to the ODE Equilibrium wall-model. A pictorial representation of the same algorithm is also provided in figure \ref{fig:algorithm_schematic}. The steps of this algorithm are described below.
	
	The wall-model internal variables and the associated ``wall variables'' stored in LES cells are initialized by the solver. As a one-time pre-processing step, the mapping routine establishes mapping between the wall-faces of the wall model and all the cells in LES grid for the rest of the simulation, as detailed in \S \ref{sec:mapping}  At each step of the simulation, the LES solver and the wall model are advanced sequentially. First, the wall model computes the wall shear stress and imposes it as a Neumann boundary condition in the LES solver at each wall face. Furthermore, the wall variables whose spatial gradients are to be computed in equation \ref{eqn:iWMx_discretized} are broadcasted to the cells in the LES grid, based on the mapping established earlier. These ``wall variables''  can now be treated as cell variables and their gradients can be readily computed using the cell-based gradient routine native to CharLES (see \S \ref{sec:surface_grad} for details). The LES is advanced in time using RK3 method and the resulting LES variables as well as the gradients of wall variables from the previous time step are passed back to the wall model, which then advances in time using the Explicit Euler method. 
	
    \begin{algorithm}[t!]
    \caption{Integral wall model LES routine}
    \label{algo:WM}
    \begin{algorithmic}[1]
    
        \State \textbf{Initialization}
        
            \State Initialize WM solver objects (assume plug flow for initializing some wall-model variables).
            \State Initialize couple objects with pointers to the LES and WM solvers.
        
        \State \textbf{Pre-processing}       
        
            \State Construct a global one-dimensional array containing wall-face information from all wall zones.        
        
            \State Using the global wall-face array, establish mapping between wall-faces of WM and cells in the LES grid. That is, based on the minimum distance criteria, assign a proper wall-face to each cell.         
        
        \State \textbf{Time integration}    
        \While{step $<$  N}  
            \State WM-to-LES BC communication: Apply $\tau_{w}$ at each wall face as Neumann BC for LES.
            \State WM-to-LES gradient data communication: Pack all wall variables whose surface gradients are to be computed into the global 1D array and broadcast to all cells in LES
            \State Compute gradients of ``wall-variables''  in each cell using cell-based gradient routine. 
            \State Advance LES one time-step.
            \State LES-to-WM BC communication: Interpolate $u_{i}$ and spatial gradient of $p$ from the LES cell-centers to the WM matching location. Spatial gradients of the integral terms ($L_{x}$, $L_{zz}$ etc.) are fed to the wall model from first off-wall cell center.
            \State Advance WM one time-step.
            \State ++step
        \EndWhile
    
    \end{algorithmic}
    \end{algorithm}

	\subsubsection{Mapping wall-faces to LES cells}\label{sec:mapping}
    We employ a routine similar to the one used in \citep{park16jcp} for establishing proper mapping of each cell in the LES grid to a wall-face. Two key differences are noted  for the current implementation; firstly, we are interested in mapping cells to wall-faces (opposed to mapping wall-model faces and LES faces in \citep{park16jcp}); secondly, we desire a one-way mapping (opposed to two-way mapping in \citep{park16jcp}) because in the current wall-model the broadcasting of data is required only from wall to the cells. It must also be emphasized that the assignment of a wall-face to a cell is not exclusive (one-to-one) i.e. each wall-face may potentially be assigned to more than one cells (one-to-many). This has a direct implication for the computation of surface gradients, as will be discussed in \S \ref{sec:surface_grad_prob}. The wallface-cell mapping routine is summarized below for completeness and shown pictorially in the \textit{pre-processing} block of figure \ref{fig:algorithm_schematic}; for details of the face mapping routine, the reader is referred to \citep{park16jcp}. 

    The LES solver first identifies the wall  zones in each processor, and writes the coordinates of wall-face centers from all these wall zones into a single global one-dimensional array. The global wall-face array is MPI-broadcasted to all the processors. For each LES cell in a given processor, the processor loops over the global array to find the wall face that is closest to that cell, based on the distance between the cell center and the wall-face center. Once the solver has looped over all the cells in all the processors, each processor will have an ordered mapping between its cells and corresponding wall faces. This mapping operation is performed as a pre-processing step and remains valid throughout the simulation, provided that the LES grid is static.

	\subsubsection{Surface gradient computation}\label{sec:surface_grad}

    As previously mentioned in \S \ref{sec:IWM_discretized}, the original formulation by Yang \emph{et al}. employs a bi-linear interpolation method to compute the surface gradients (i.e., gradient of quantities defined on the wall). However, the bi-linear interpolation is not feasible in an unstructured finite-volume solver. Specifically, in the CharLES code, pre-built structures for the wall-face connectivity are not available, and therefore accessing the neighboring wall-faces is not a trivial task. 

\begin{figure}[t!]
		\begin{minipage}{0.45\textwidth} 
			\includegraphics[width=\linewidth]{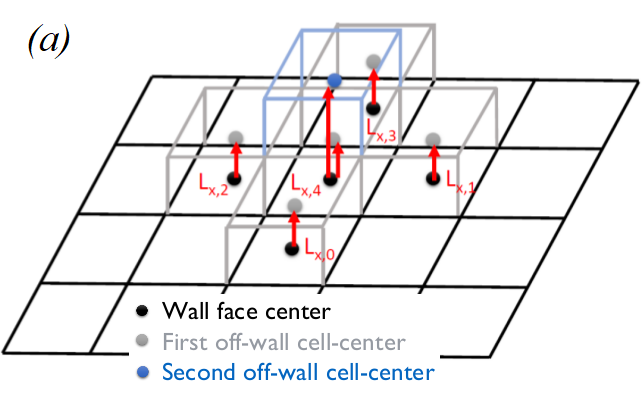}    
		\end{minipage}    
		\hspace{\fill}  
		\begin{minipage}{0.45\textwidth} 
			\includegraphics[width=\linewidth]{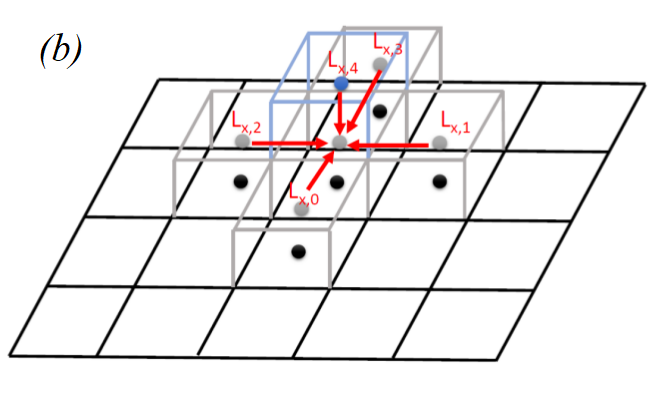}    
		\end{minipage}    
		
		\vspace{0.75cm}
		\begin{minipage}{0.45\textwidth} 
			\includegraphics[width=\linewidth]{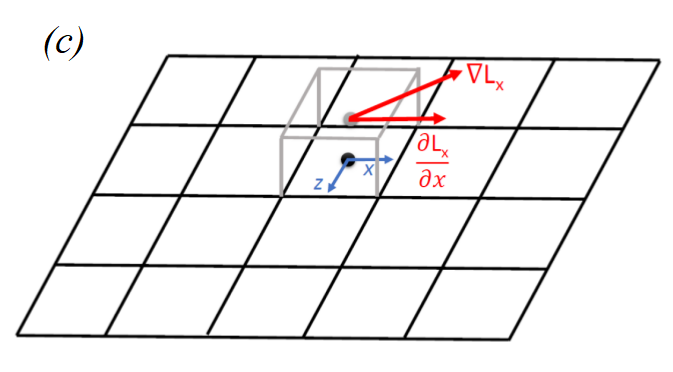}    
		\end{minipage}    
		\hspace{\fill} 
		\begin{minipage}{0.45\textwidth} 
			\includegraphics[width=\linewidth]{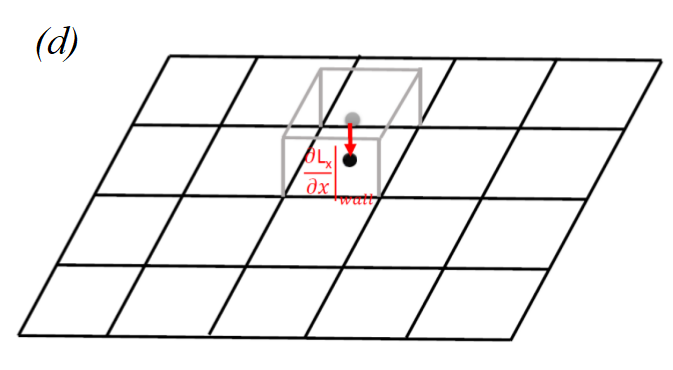}    
		\end{minipage}    
		
        \caption[]
        {\small Surface-gradient computation for wall variables.  (a) Step 1: CVs are filled with wall variables from the closest wall face, (b) Step 2: Gradient of the integral terms are computed in CVs, using the cell-based gradient routine native to the flow solver, (c) Step 3: Derivative along the desired local coordinate is obtained by taking the projection of $\nabla L_{x}$ onto that direction, (d) Step 4: $\frac{dL_{x}}{dx}$ is copied from the first off-wall cell-center to the wall-face center. Note that the axes $x$ and $z$ define the local wall coordinate system.} 
        
        \label{fig:surface_gradient}
	\end{figure}

    CharLES has a built-in routine for computing cell gradients, where data from the first-neighbor cells is used (Green-Gauss second-order reconstruction). As an alternative to the bi-linear interpolation method for computing surface gradients, the wall-face quantities (integral terms) are mapped to the cell-centers of control volumes (CVs) based on the wallface-cell mapping established in the pre-processing step, and their gradients are computed using the cell-gradient routine. Subsequently, the gradients from the first off-wall cell are copied back to the wall-faces, where they can be used by the wall-model. It should be noted that in our current implementation each of the integral terms in 3D ($L_{x}$, $L_{xx}$, $L_{z}$, $L_{zz}$, $L_{xz}$) is treated as a scalar quantity (five scalar integral quantities per wall face) and these are copied to the cell-centers as scalars. The gradients are then computed based on these scalars. For instance, to obtain $\frac{dL_{x}}{dx}$, we compute $\nabla L_{x}$ (where $L_{x}$ is a scalar) from values of $L_{x}$ in neighboring cells and then take the component of this gradient vector  $\nabla L_{x}$ in the desired local $x$-direction. This process is shown step by step in figure \ref{fig:surface_gradient}. Similarly gradients are computed separately for rest of the 4 integral terms and their derivatives in the desired direction are obtained.

	\subsection{Challenges associated with surface-gradient computation in unstructured solver}\label{sec:surface_grad_prob}
	
	As alluded to previously, the computation of surface gradients in an unstructured solver presents a couple of challenges. The first (and more critical) problem pertains to the gradient of integral terms in the 3D  integral momentum equation set given below, which is the generalization of equation \ref{eqn:iWMx}. 
    
    \begin{equation}\label{eqn:iWMx_3D}
    \begin{split}
        \frac{\partial L_{x}} {\partial t} + \frac{\partial L_{xx}}{\partial x} + \frac{\partial L_{xz}}{\partial z}
        -U_{LES} \left(\frac{\partial L_{x}}{\partial x} + \frac{\partial L_{z}}{\partial z}\right) & = 
        \frac{1}{\rho} \left[-\frac{\partial p}{\partial x} h_{wm} +\tau_{h_{wm},x} - \tau_{w,x}\right] \\
        \frac{\partial L_{z}} {\partial t} + \frac{\partial L_{xz}}{\partial x} + \frac{\partial L_{zz}}{\partial z}
        -W_{LES} \left(\frac{\partial L_{x}}{\partial x} + \frac{\partial L_{z}}{\partial z}\right) & = 
        \frac{1}{\rho} \left[-\frac{\partial p}{\partial z} h_{wm} +\tau_{h_{wm},z} - \tau_{w,z}\right]
    \end{split}
    \end{equation}

    \noindent where,
    \begin{equation}\label{eqn:Lx_3D}
    L_{x} \equiv \int_{0}^{h_{wm}} u dy  , \; \; \;  
    L_{xx}\equiv \int_{0}^{h_{wm}} u^2 dy, \; \; \;
    L_{z} \equiv \int_{0}^{h_{wm}} w dy  , \; \; \;  
    L_{zz}\equiv \int_{0}^{h_{wm}} w^2 dy, \; \; \;
    L_{xz}\equiv \int_{0}^{h_{wm}} u w dy .
    \end{equation}

		
        

	\begin{figure}
			\includegraphics[trim=0 0 0 30,clip,width=\linewidth]{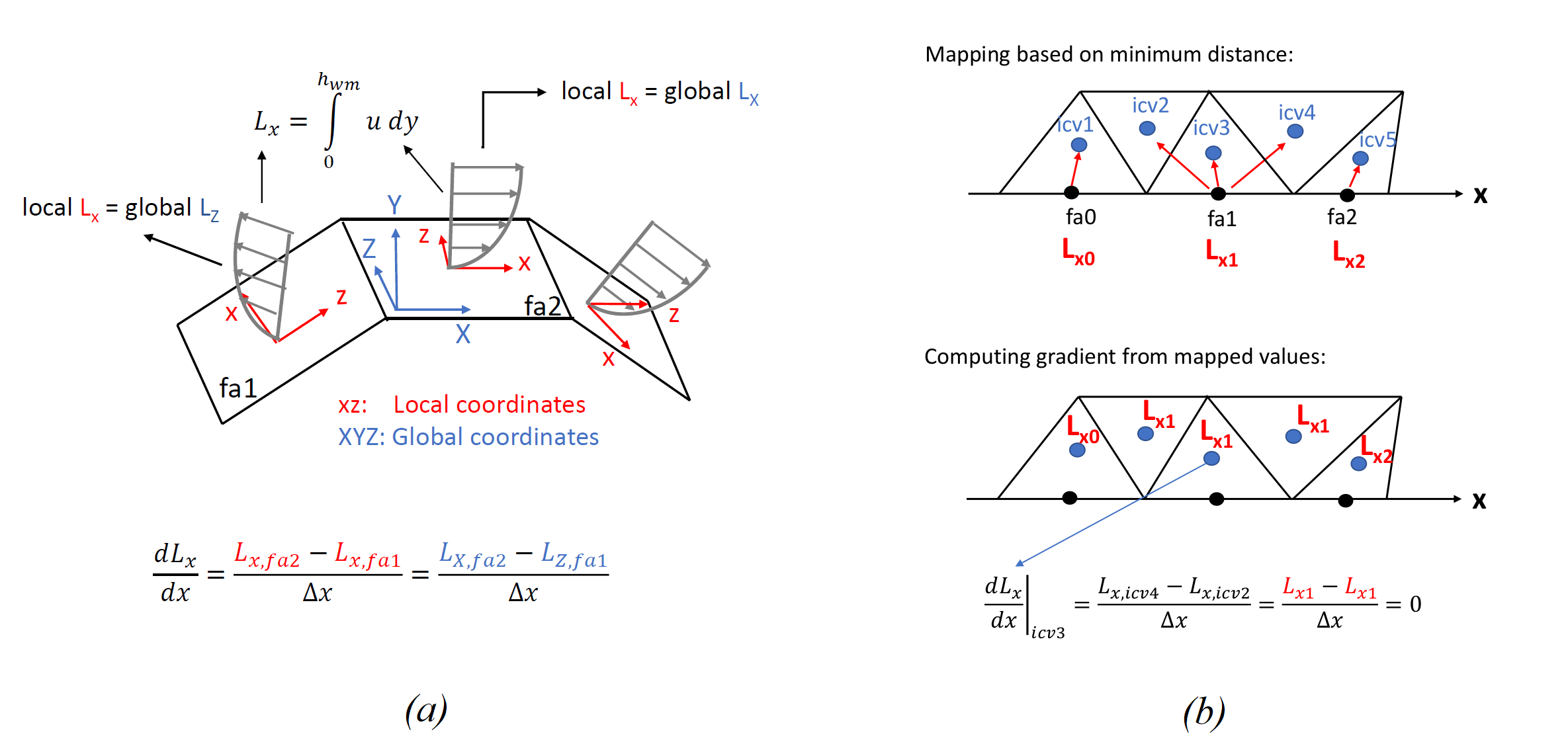} 
		
        \caption[]
        {\small Challenges associated with surface gradient, (a) Schematic of the surface gradient problem related to drastic change of the local wall coordinate system, (b) Non-exclusive (one-to-many) mapping of wall-variables to cells and the resulting inaccurate computation of gradients by the cell-based gradient routine.} 
        
        \label{fig:surface_gradient_prob}
	\end{figure}

    \noindent It should be reiterated that $x$, $z$ here are the \emph{local} wall-parallel coordinates which in principle can be chosen arbitrarily, and that $y$ is aligned with the local wall-normal direction. Specifically, the transformation of gradients of integral terms to a different coordinate system poses a challenge and may result in inaccurate results if not done using a proper tensor transformation rule. The algorithm discussed in \S \ref{sec:surface_grad} works well under the ideal conditions of a regular Cartesian grid with small wall-surface curvature. However, there is a key limitation associated with this approach; it assumes that the local $x$ and $z$ axes remain fairly constant in the neighborhood of a given wall-face. This assumption does not always hold, especially in regions with rapid change in the wall curvature or at junctures. For instance, computing $\nabla L_{x}$ with a drastically changing local $x$-direction in adjacent wall faces would be akin to computing this gradient based on different scalar quantities from each wall-face. This point is illustrated through a pathological scenario in figure \ref{fig:surface_gradient_prob} $(a)$, where the local $x$-axis rotates through a $90^{0}$ angle in the wall parallel plane, in going from \textit{face 1} to \textit{face 2}. Here the global coordinate axes $X$-$Z$ are respectively aligned with the local $x$ and $z$ axes of \textit{face 2}. Consequently, the integral term $L_{x}$ in local coordinates for \textit{face 1} corresponds to $L_{Z}$ in the global coordinate system.  Therefore, computation of $\nabla L_{x}$ by the cell-based gradient routine (which operates in the global coordinate system) involves the $L_{Z}$ term from \textit{face1} and $L_{X}$ term from \textit{face2}. This is where the current implementation is prone to error. As a quick workaround, we propose fixing the local coordinates for all the wall faces in the pre-processing step, such that any abrupt spatial changes in local coordinate axes are suppressed.


    
    The second problem with the cell-based gradient approach is that it would give accurate values for gradients only for a hex-dominant grid near the wall, like that seen in figure \ref{fig:surface_gradient}, where association between wall faces and first off-wall cells is unique based on the minimum-distance criteria. However, for a grid with tetrahedral elements next to the wall, the mapping of variables from the wall-face to the cells would not be unique for the first layer of cells adjacent to the wall, resulting in inaccurate gradient computation. An illustrative example of this scenario is shown in figure \ref{fig:surface_gradient_prob}$(b)$ for a 2D grid with triangular elements. Here, the gradient in icv3 is erroneously calculated to be zero because all the neighboring cells of icv3 get values from the same face fa1 based on the minimum-distance criteria. A quick fix for this problem is to use a spatial filter for the computed wall-gradients in the wall-parallel plane, before feeding them back to the wall-faces.

	\subsection{Predetermination of Quadrature points and associated error}	
	
\begin{figure}[t]
    	\centering
    	\begin{subfigure}[b]{0.44\textwidth}
    		\includegraphics[trim=100 230 100     230,clip,width=\textwidth]{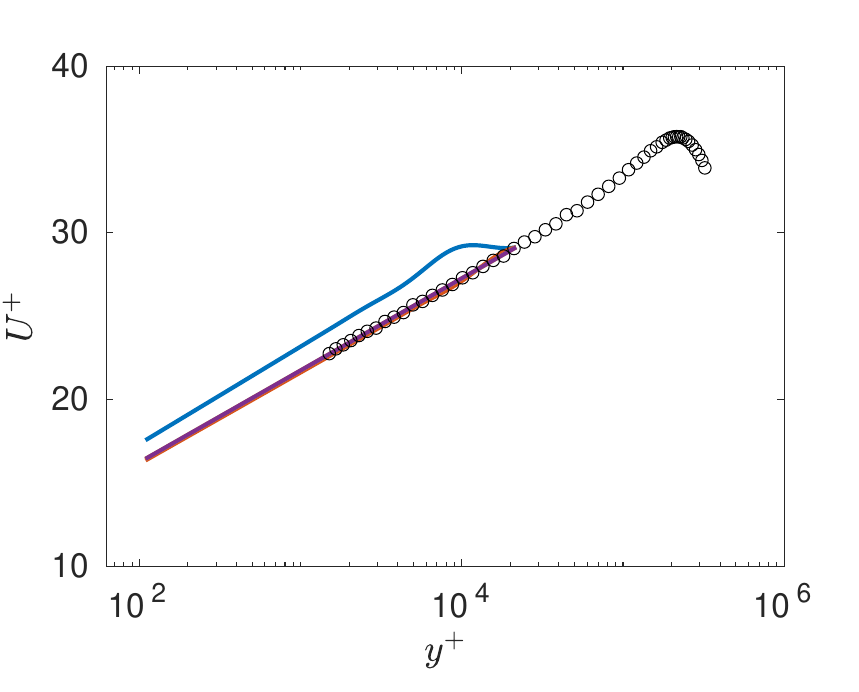}
    	    \caption{}
    	\end{subfigure}
    	~
    	\begin{subfigure}[b]{0.44\textwidth}
    		\includegraphics[trim=100 230 100     230,clip,width=\textwidth]{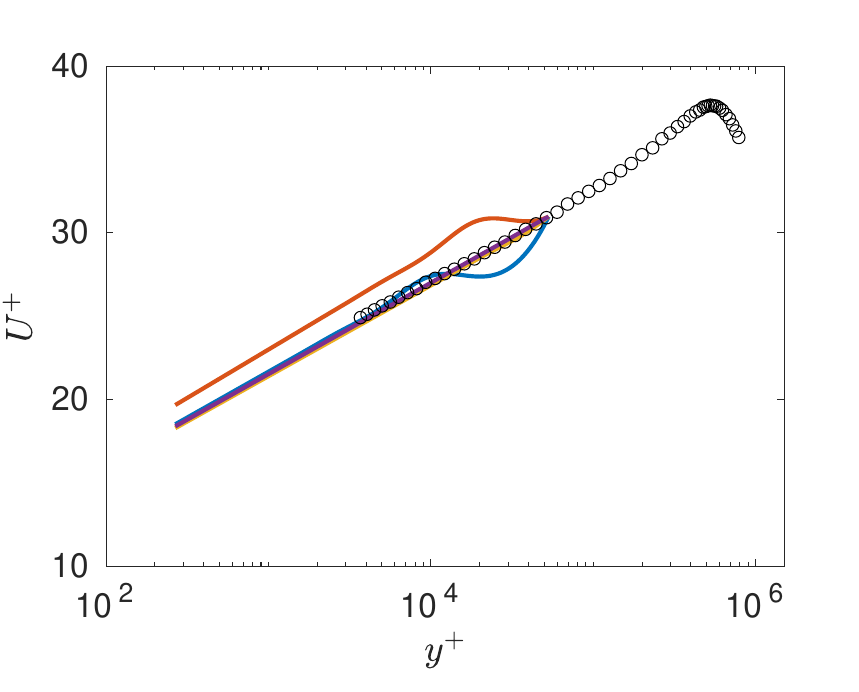}
    	    \caption{}
    	\end{subfigure}

    	\caption{\small Spurious predicted velocity profile for Gauss-quadrature based method resulting from insufficient number of quadrature points (a) $Re_{\tau} = 210,000 $, (b) $Re_{\tau} = 525,000 $. The colors represent different number of quadrature points used in the wall model; blue, $n=70$; red, $n=100$; yellow, $n=200$; violet, $n=300$. The circles are from the turbulent pipe flow experiment \citep{zagarola98}. }
    	\label{fig:GQ_spurious_profile}
    \end{figure}

	An important caveat of using Gauss-quadrature technique for ODE wall-model implementation is that the user needs to prescribe the number of quadrature points for integral evaluation beforehand. Although a similar requirement exists for the FV-based method in terms of choosing the number of wall-model grid points, however it must be emphasized that the Gauss-quadrature being a spectral method utilizes basis functions to approximate the integrand, and as such is more prone to producing a spurious velocity profile if insufficient quadrature points are used to approximate the underlying velocity profile. This is especially the case for extremely high Reynolds number, where the number of quadrature point requirement increases significantly, as depicted in figure \ref{fig:GQ_spurious_profile}.

	One possible remedy is to compute the $L_{2}$ error of the predicted velocity profile against an extremely fine ($n>200$) case which acts as a reference velocity profile. However, this additional check imposes a significant overhead and might possibly degrade the performance of this method. It should however be noted that the cost of Gauss-quadrature is significantly lower than FV method even for $n_{GQ}>>n_{FV}$ due to the cost scaling shown in figure \ref{fig:cost_vs_n}, allowing us to prescribe a large number of quadrature points and still get better performance than the traditional FV method. This point will be further expanded upon in \S\ref{sec:GQ_perform}.

	\section{Performance considerations and validation}\label{sec:performance_validation}
	


    

	\subsection{Validation of wall models}

\begin{figure}[t]
    	\centering
    	\begin{subfigure}[b]{0.46\textwidth}
    		\includegraphics[trim=80 200 100 250,clip,width=\textwidth]{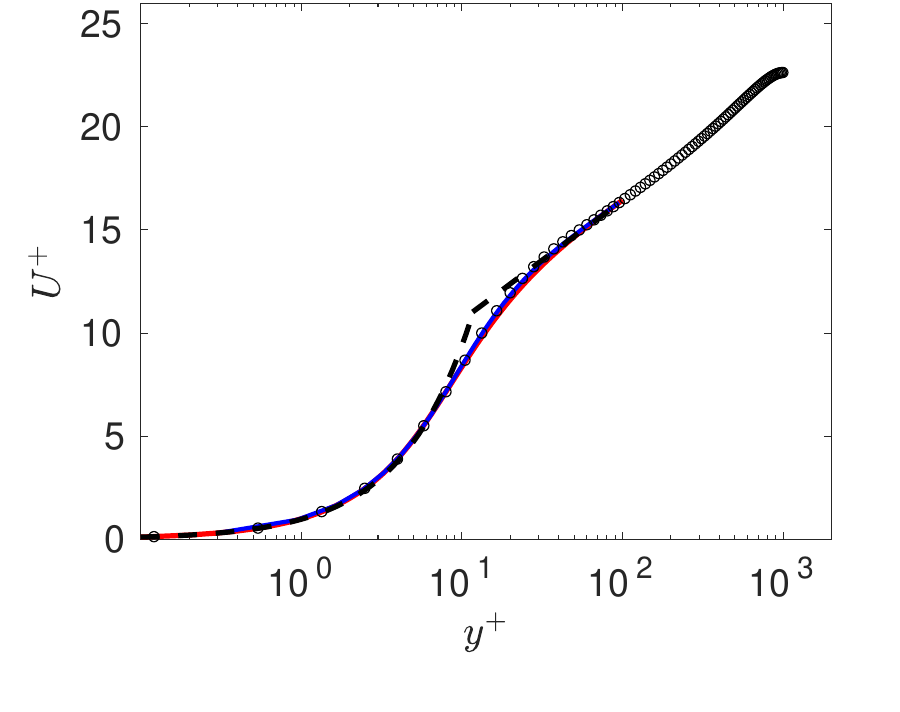}
    		\caption{}
    		\label{fig:apriori_meanprofile}
    	\end{subfigure}
    	~
    	\begin{subfigure}[b]{0.46\textwidth}
    		\includegraphics[trim=80 200 100 250,clip,width=\textwidth]{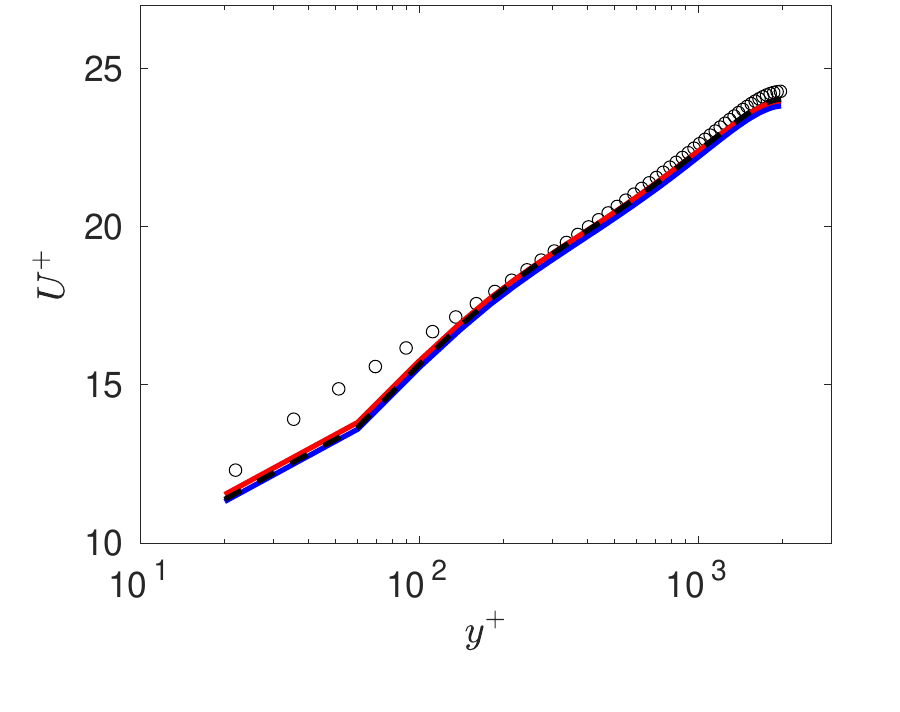}
    		\caption{}
    		\label{fig:aposteriori_meanprofile}
    	\end{subfigure}

    	\caption{\small Validation of the wall models; (a)  \textit{a priori} mean velocity profile from the wall models; (b) \textit{a posteriori} mean velocity profiles from the wall-modeled LES solution. Blue solid line, finite-volume based ODE WM; red solid line, Gauss-Quadrature based ODE WM; black dashed line, integral WM; circles in (a), channel flow DNS \citep{graham16} at $Re_{\tau}=1000$; circles in (b), channel flow DNS \citep{hoyas06} at $Re_{\tau}=2000$. }
    	\label{fig:validation}
    \end{figure}

    Both the wall-models (GQWM and integral NEQWM) were validated in \textit{a priori} as well as \textit{a posteriori} settings. For \textit{a priori} test, DNS of turbulent channel flow at $Re_{\tau} = 1000$ by Graham \textit{et al.} \citep{graham16} was used as reference, due to its ease of accessibility and suitability for testing integral NEQWM; whereas for \textit{a posteriori} test, WMLES was performed at $Re_{\tau} = 2000$, where $Re_{\tau}$ is the Reynolds number based on the channel half-height $\delta$, friction velocity $u_{\tau}$, and kinematic viscosity $\nu$. The size of the computational domain in WMLES is $L_{x}=25\delta$, $L_{y}=2\delta$, and $L_{z}=10\delta$, where $x$, $y$, and $z$ represent the streamwise, wall-normal, and spanwise directions, respectively. The LES grid contains a total of $2.7$ million grid points with $(N_{x},N_{y},N_{z}) = (180,100,150)$, resulting in a uniform grid spacing of $(\Delta x^{+},\Delta y^{+},\Delta z^{+}) = (277.8,40,133.3)$ in wall units. The flow configuration is comparable to the WMLES in \citep{park16prf} and \citep{park14pof}. The reference used for the \textit{a posteriori} test is the incompressible channel DNS of Hoyas and Jiménez \citep{hoyas06}. Figures \ref{fig:apriori_meanprofile} and \ref{fig:aposteriori_meanprofile} respectively show the mean velocity profiles from these \textit{a priori} and \textit{a posteriori} tests for the two wall models, and their comparison with FV-based ODE EQWM and DNS. We observe almost perfect agreement of the two wall models with the FV-based wall model and a reasonable agreement with the channel flow DNS. \textit{A priori} tests for GQWM were performed for a wide range of $Re_{\tau}$ ($\sim 10^{3} - 10^{6}$) based on various canonical cases available in the literature, namely turbulent pipe flow experiment \citep{zagarola98}, turbulent boundary layer experiment \citep{osterlund99}, and turbulent channel flow DNS \citep{lee15}. This allowed us to analyze the cost scaling of GQWM with $Re_{\tau}$ and the number of quadrature points $n$ which will be discussed in \S \ref{sec:GQ_perform_apriori}. The matching location for the wall model in all of these validation tests was set to $10\%$ of the boundary layer thickness i.e. $h_{wm}=0.1\delta$. The MATLAB implementation for the \textit{a priori} validation in channel flow for all three wall models can be found in the authors' GitHub repository \url{https://www.github.com/imranhayat29/Wall-Models-for-LES}.

	\subsection{Performance characteristics of Gauss Quadrature method}\label{sec:GQ_perform}

	A quadrature based method requires significantly lesser number of quadrature points than the number of cells required by the FV method to approximate a polynomial of the given order. Furthermore, for a given number of quadrature points or cells $n$, the FV method involves inversion of size $n \times n$ matrix whereas the quadrature method involves only $n$ function evaluations. So, even if the number of quadrature points for GQ and cells for FV are comparable number, the GQ method is expected to outperform FV method in terms of the computational cost of the wall model. In the following subsections, we evaluate the performance characteristics of these methods in  \textit{a priori} and \textit{a posteriori} settings. Note that \textit{a priori} tests are quite insightful in assessing the cost difference between the two methods, in particular identifying the cause of that difference and its scaling with $Re_{\tau}$ and $n$. \textit{A posteriori} test, on the hand, sheds light on how the two methods perform in parallelized settings, in particular how the load-balancing among the processors (or lack thereof) for the wall-model solver affects the two methods.

    \subsubsection{\textit{A priori} performance}\label{sec:GQ_perform_apriori}

\begin{figure}[t]
    	\centering
    	\begin{subfigure}[b]{0.49\textwidth}
    		\includegraphics[trim=50 200 50 200,clip,width=\textwidth]{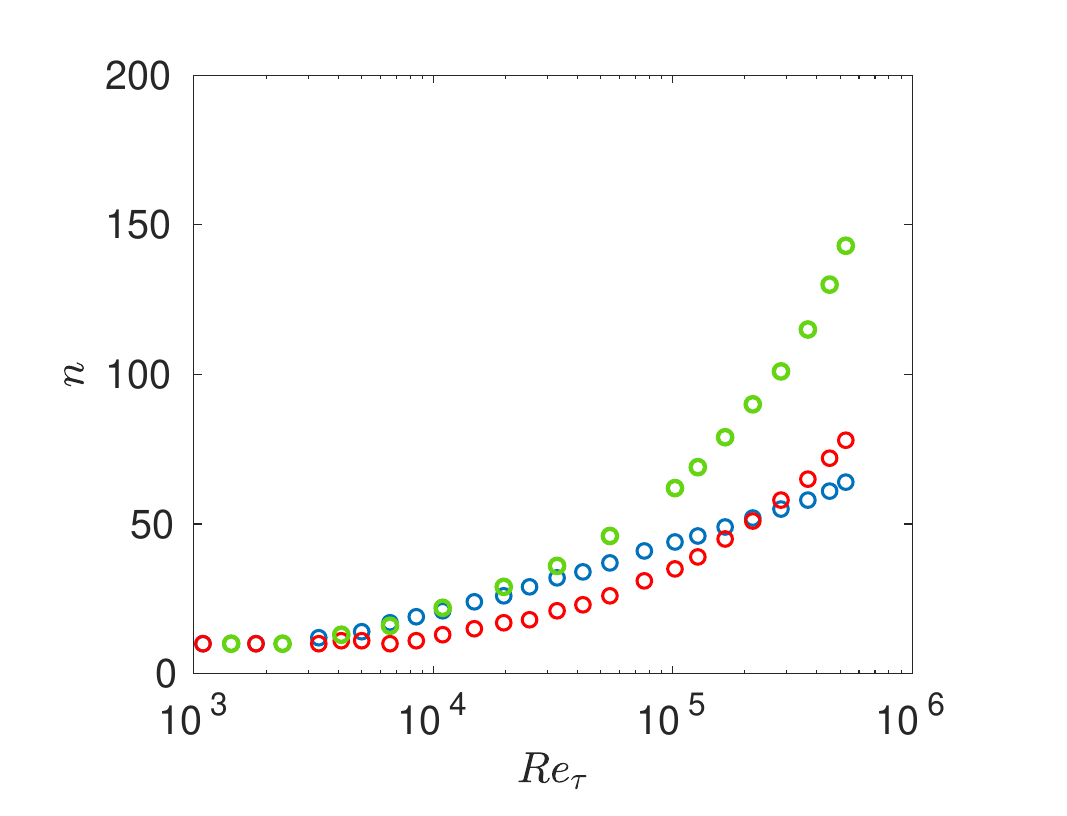}
    		\caption{}
    		\label{fig:optncv_tpf}
    	\end{subfigure}
    	~
    	\begin{subfigure}[b]{0.49\textwidth}
    		\includegraphics[trim=50 200 50 200,clip,width=\textwidth]{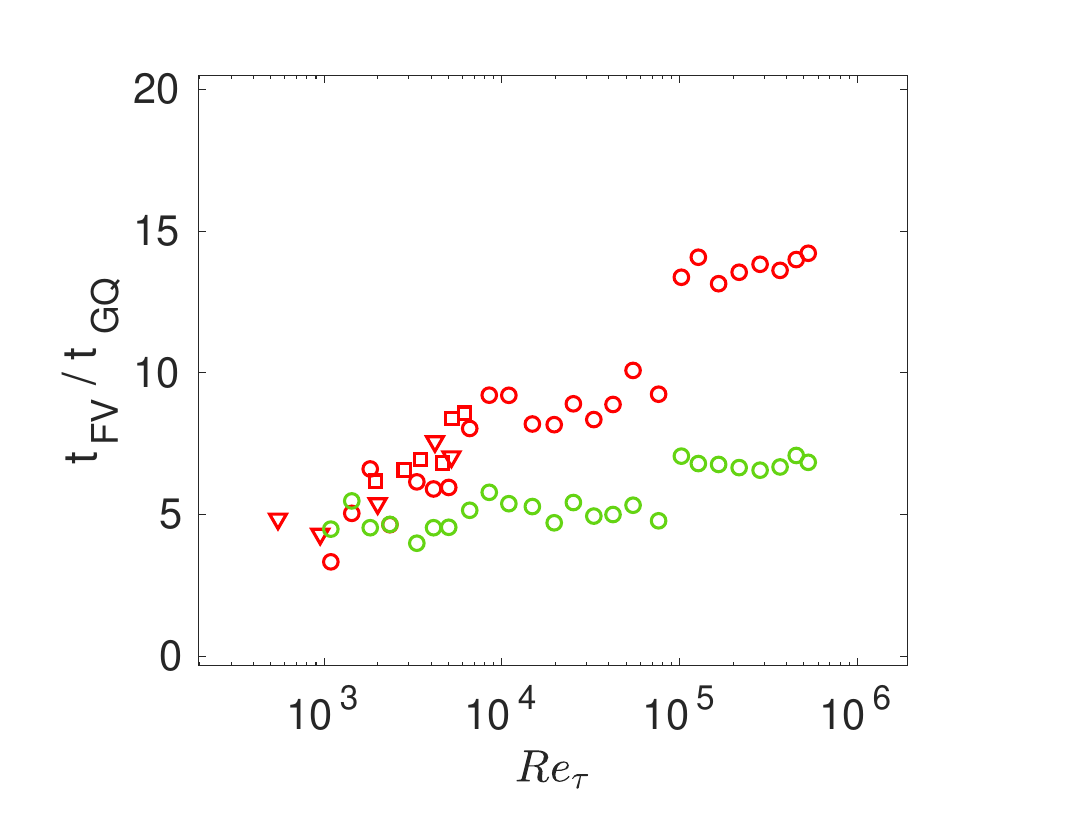}
    		\caption{}
    		\label{fig:speedup_tpf}
    	\end{subfigure}

    	\caption{\small (a) Optimal number of points required to achieve less than $3\%$ error in the predicted $\tau_{w}$ at different $Re_{\tau}$; (b) Speedup for GQ over FV approach. Symbols in (b) are from various studies: triangles, channel flow DNS \citep{lee15}; squares, turbulent boundary layer experiment \citep{osterlund99}; circles, turbulent pipe flow experiment \citep{zagarola98}. Colors represent different methods: blue, FV; green, GQ with linear transformation; red, GQ with non-linear transformation (clustering).}
    	\label{fig:apriori_tpf}
    \end{figure}

    Figure \ref{fig:optncv_tpf} shows the scaling of the optimal number of points $n$ in the wall model with friction Reynolds number, where optimal $n$ is defined as the $n$ which produces an error of $\leq 3\%$ in the predicted wall shear stress magnitude compared with the reference $\tau_{w}$ from DNS. For GQ based model with linear transformation given by equation \ref{eqn:lin_transf}, the optimal number of quadrature points $n$ scales strongly with Reynolds number as $Re_{\tau}^{0.5}$ compared with $Re_{\tau}^{0.2}$ scaling for the FV based model. However, the non-linear transformation in equation \ref{eqn:nonlin_transf} with clustering of quadrature points close to the wall results in a decrease in $n$ required to achieve the error tolerance in wall stress, and a weaker scaling with $Re_{\tau}$, as seen in figure \ref{fig:optncv_tpf}. Figure \ref{fig:speedup_tpf} shows the corresponding speedup of GQ over FV based model, for both linear and clustered transformations. The speedup here is defined as the ratio of time taken by the FV model to that by the GQ model for a given Reynolds number. Note that the time here is the average wall-clock time recorded over several hundred runs of each wall model. It is evident that even with the linear transformation which requires a large number of quadrature points especially at very high Reynolds number, the GQ method shows a significant speedup ($~4-6$ times) over FV method. With the clustered GQ, the speedup is even higher, reaching up to $14$ times for the highest Reynolds number considered in this \textit{a priori} study. Furthermore, the speedup appears to be a strong function of the Reynolds number when clustered transformation is used, whereas for linear transformation the speedup remains fairly constant with increase in Reynolds number. From this point onward, we will only analyze the clustered GQ method.

\begin{figure}[t]
    	\centering
    	\begin{subfigure}[b]{0.49\textwidth}
    		\includegraphics[trim=50 200 50 200,clip,width=\textwidth]{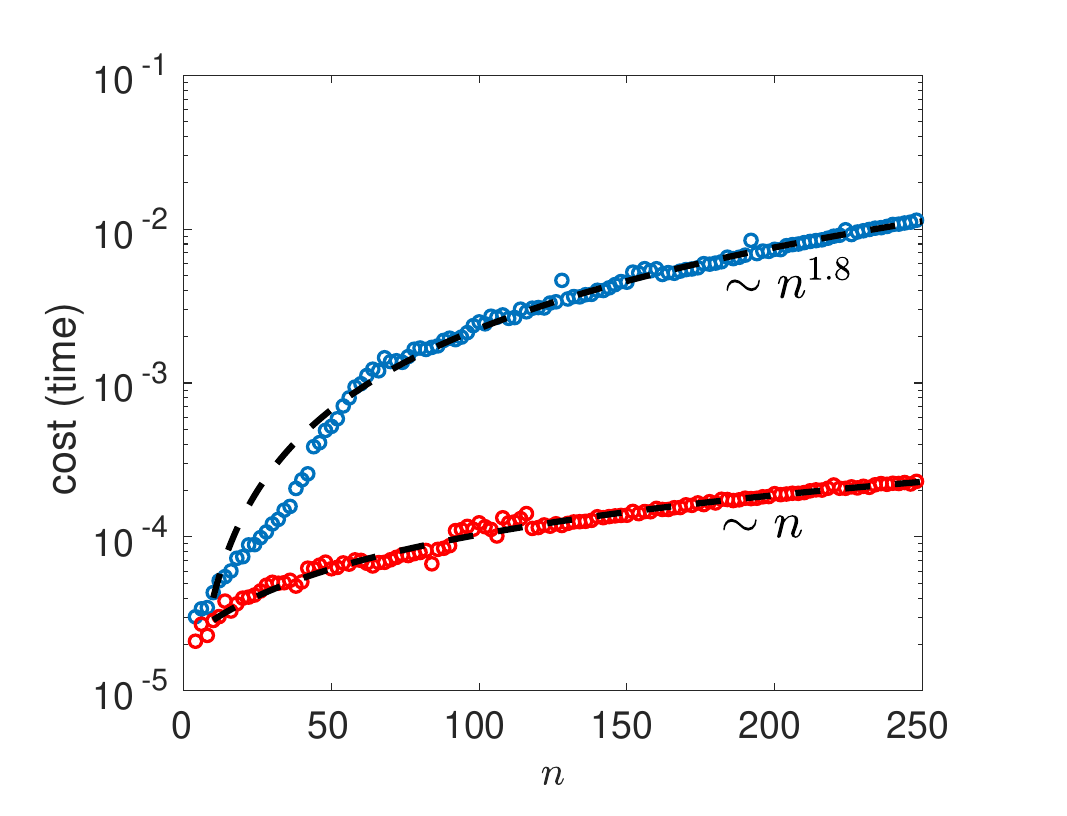}
    		\caption{}
    		\label{fig:cost_vs_n}
    	\end{subfigure}
    	~
    	\begin{subfigure}[b]{0.49\textwidth}
    		\includegraphics[trim=50 200 50 200,clip,width=\textwidth]{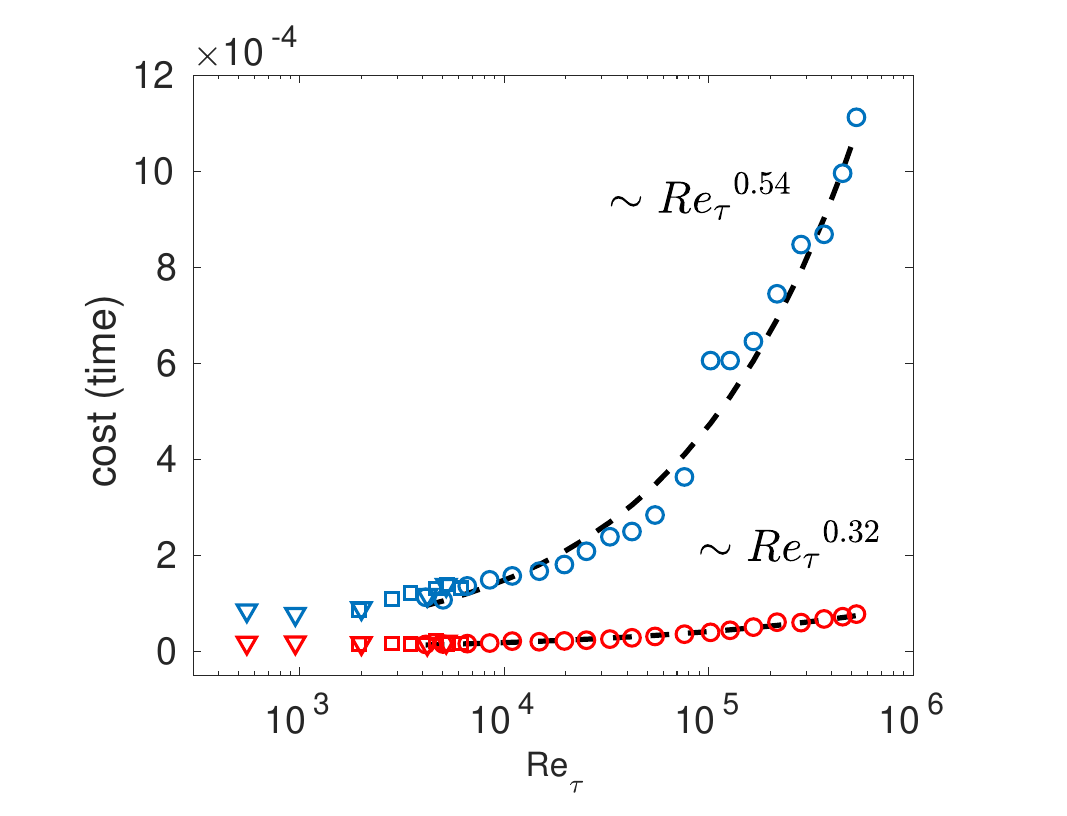}
    		\caption{}
    		\label{fig:cost_vs_Retau}
    	\end{subfigure}

    	\caption{\small \textit{A priori} cost analysis; (a) cost scaling with the number of points used in the wall model; (b) cost scaling with friction Reynolds number. Colors in both (a) and (b), and symbols in (b) are the same as those described in figure \ref{fig:speedup_tpf}. }
    	\label{fig:apriori_cost_tpf}
    \end{figure}

    To understand why GQ method shows speedup over FV method, we need to look at the cost scaling of the two methods with $n$. Figure \ref{fig:cost_vs_n} shows this cost scaling for the clustered GQ and FV wall models. It is observed that the computational cost of FV method has almost an $O(n^2)$ scaling whereas GQ method scales as $O(n)$. This is consistent with the fact stated earlier that the FV method involves inversion of size $n \times n$ matrix whereas the quadrature method involves $n$ function evaluations. This cost scaling explains the increase in speedup for GQ at higher $Re_{\tau}$ in figure \ref{fig:speedup_tpf}; as $Re_{\tau}$ is increased, more points $n$ are required in the wall model to capture the near-wall region in $y^{+}$ units, correspondingly the scaling in figure \ref{fig:cost_vs_n} at higher $n$ comes into play, thus widening the cost gap between the two methods. Figure \ref{fig:cost_vs_Retau} illustrates this widening gap for increasing $Re_{\tau}$, where the cost of FV method scales as $Re_{\tau}^{0.54}$ and that of GQ method scales as $Re_{\tau}^{0.3}$. Considering that in a real simulation, the wall modeling calculations are done for a large number of wall faces and at each time step of the simulation, this difference in wall-modeling cost between the two models, especially for high Reynolds number flows, can result in a significant computational cost reduction for the GQ based wall model.

    \subsubsection{\textit{A posteriori} performance}\label{sec:GQ_perform_aposteriori}

\begin{figure}[t]
    	\centering
    	\begin{subfigure}[b]{0.46\textwidth}
    		\includegraphics[trim=40 180 50 220,clip,width=\textwidth]{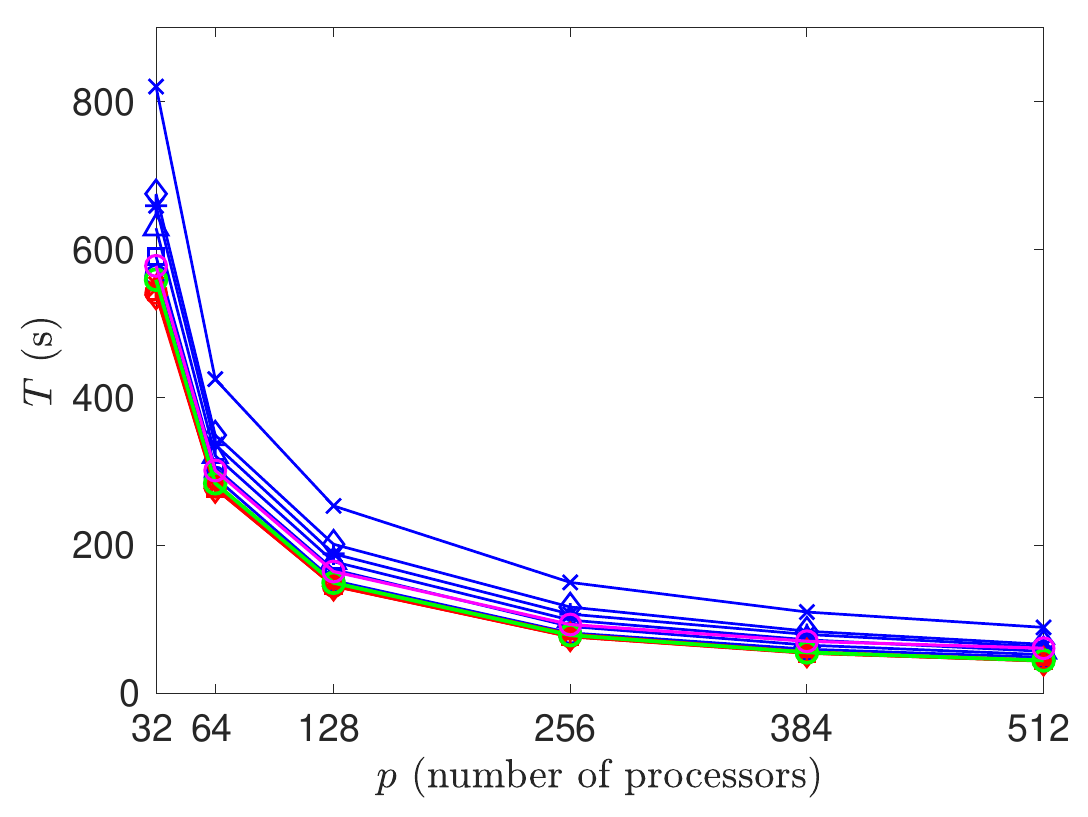}
    		\caption{}    		
    	\label{fig:aposteriori_WMcost}	
    	\end{subfigure}
    	~  
        \begin{subfigure}[b]{0.46\textwidth}    	    \includegraphics[trim=40 180 50     220,clip,width=\textwidth]{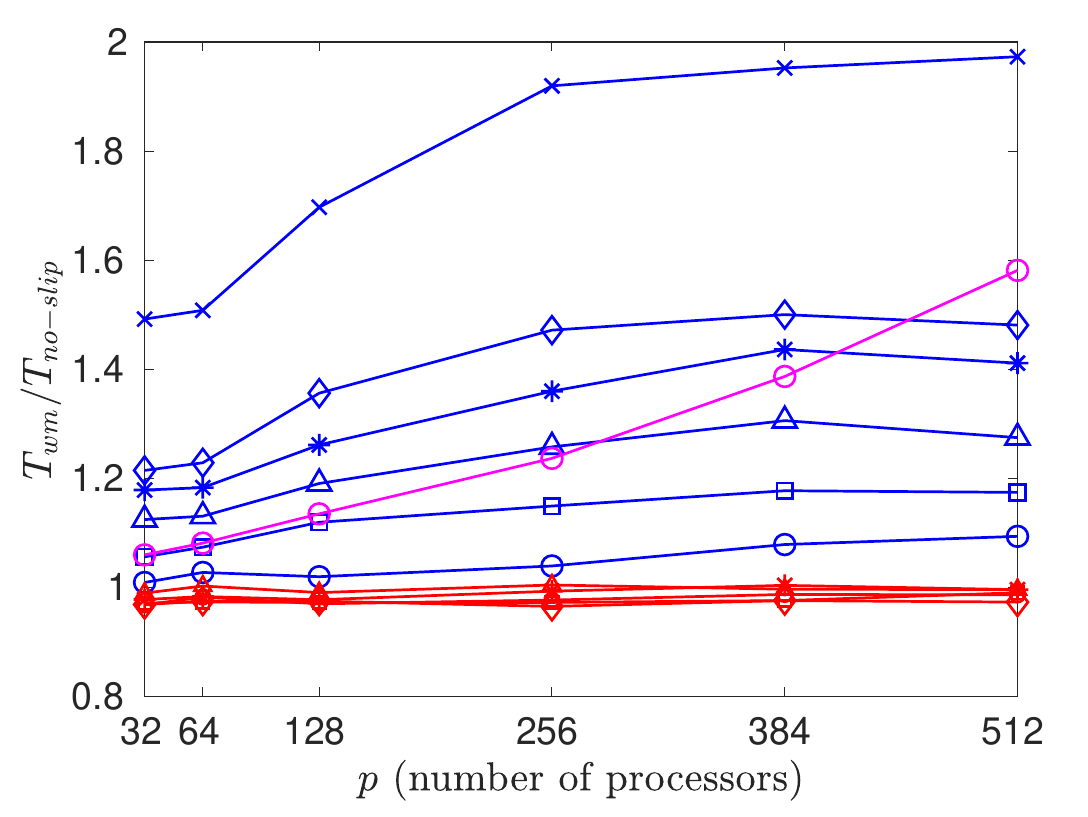}
    		\caption{}        
    	\label{fig:aposteriori_add_WMcost}
    	\end{subfigure}
    	
    	\caption{\small  (a) \textit{A posteriori} computational cost of simulation; (b) additional cost of the wall model compared to the no-slip simulation. Colors represent the type of wall model used in LES: blue, FV ODEWM; red, GQ ODEWM; magenta, integral NEQWM; green, no-slip. Symbols represent the number of quadrature or grid points used in the ODE wall models: circle, $n=10$; square, $n=20$; triangle, $n=30$; asterisk, $n=40$; diamond, $n=50$; cross, $n=100$. (Note that $n$ is meaningless for integral NEQWM and no-slip.)}    	
    	\label{fig:aposteriori_cost_GQ}
    \end{figure}

    Figure \ref{fig:aposteriori_cost_GQ} shows \textit{a posteriori} computational cost analysis for all the wall models considered in this study. Figure \ref{fig:aposteriori_WMcost} compares the simulation time of LES for GQ and FV based wall models with respect to the stand-alone LES without a wall model (hereinafter referred to as no-slip LES). Note that the simulation time here is defined as the average wall-clock time required to complete 1000 time-steps of the simulation. As expected, the simulation time decreases with increase in number of processors. More importantly, the simulation time for FVWM is substantially higher than no-slip simulation, whereas that for GQWM is identical to no-slip simulation. To make this distinction more clear, we plot the ratio $T_{wm}/T_{no-slip}$ in figure \ref{fig:aposteriori_add_WMcost}. This essentially provides information on the \textit{additional cost} incurred by using a wall model in the LES compared to stand alone LES without a wall model. It is clear from the figure that the GQWM has virtually no additional cost over the no-slip simulations, irrespective of the number of quadrature points $n$ used in the wall model or the number of processors $p$ used in the simulation. In contrast, the FVWM incurs a significant cost overhead, which is increases with both $n$ and $p$. To put this in perspective, for the highest number of points ($n = 100$) and the highest number of processors ($p = 512$) considered in this study, an additional wall-model cost of up to $100\%$ of the no-slip simulation was observed for the FV method.

\begin{figure}[t]
    	\centering
    	\begin{subfigure}[b]{0.46\textwidth}
    		\includegraphics[trim=50 150 50 160,clip,width=\textwidth]{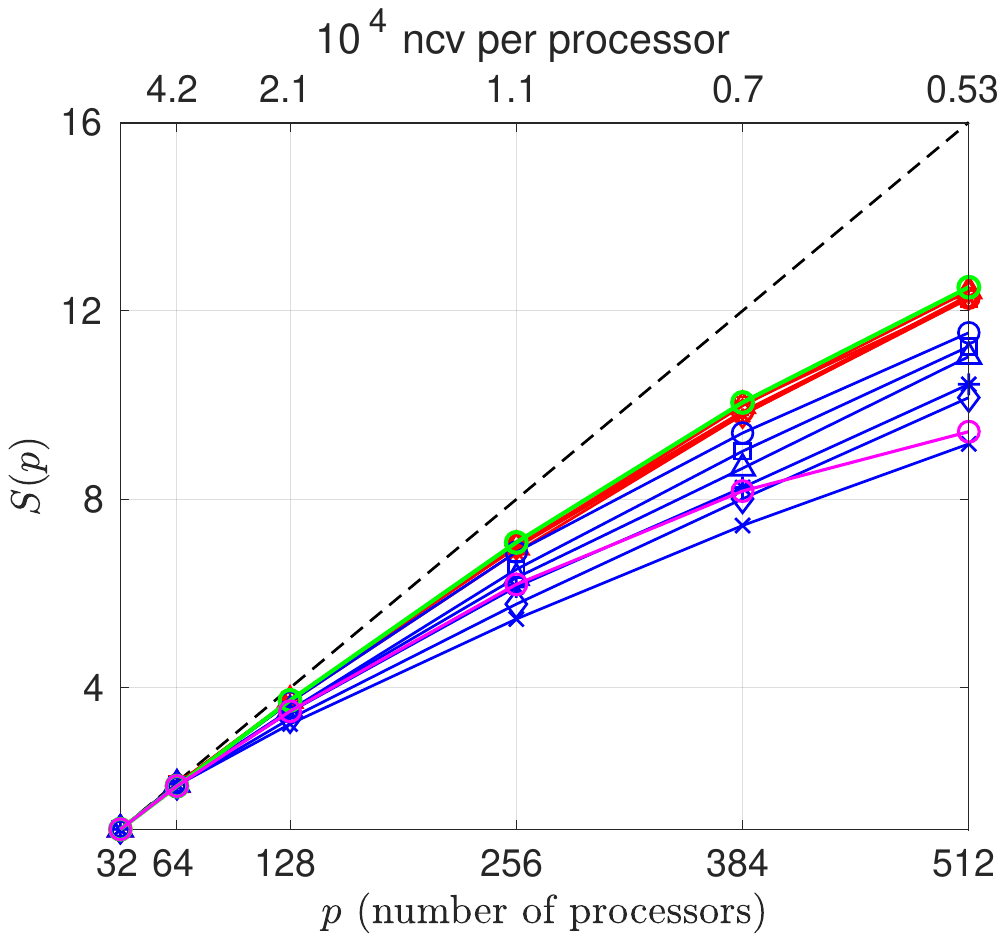}
    		\caption{}
    		\label{fig:speedup_GQ_aposteriori}
    	\end{subfigure}
    	~
    	\begin{subfigure}[b]{0.46\textwidth}
    		\includegraphics[trim=50 150 50 160,clip,width=\textwidth]{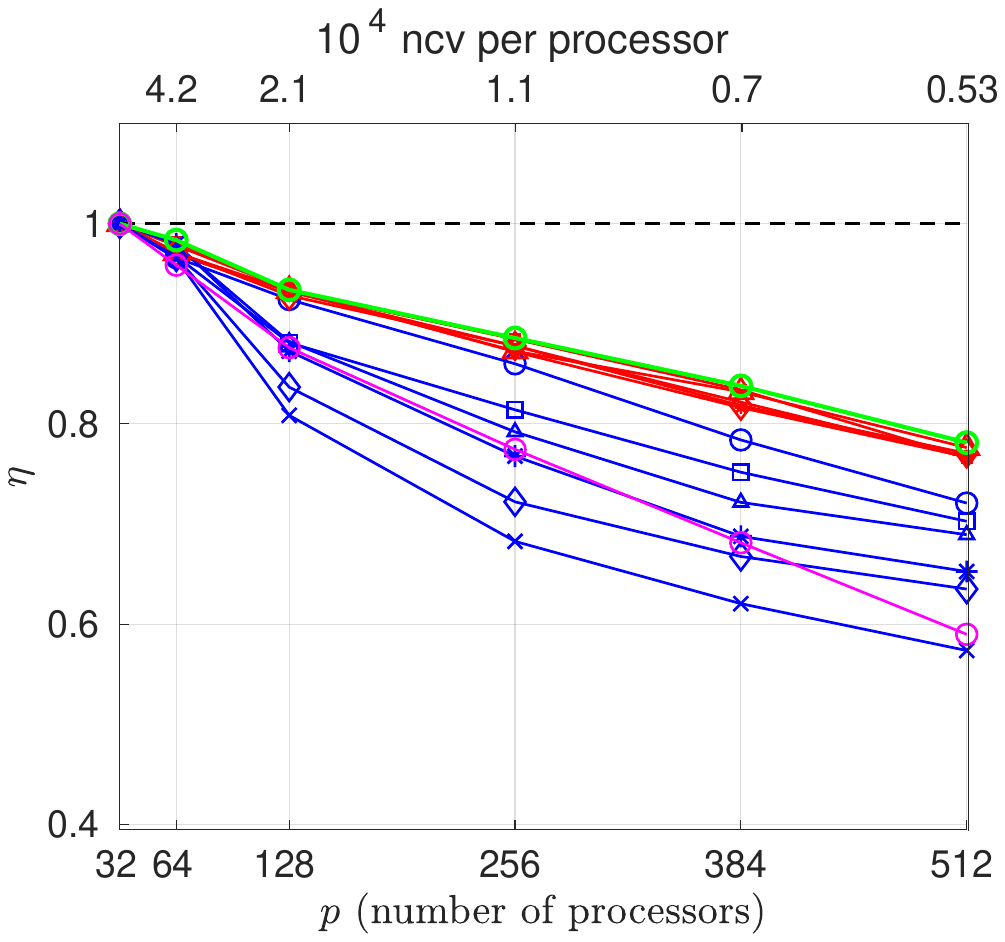}
    		\caption{}
    		\label{fig:pareff_GQ_aposteriori}
    	\end{subfigure}
    	
    	\caption{\small  \textit{A posteriori} parallel performance analysis of GQ vs FV based ODE EQWM;  (a) speedup; (b) parallel efficiency. Dashed line shows the ideal (linear) scaling; symbols and colors are the same as described in figure \ref{fig:aposteriori_cost_GQ}.}
    	\label{fig:aposteriori_par_GQ}
    \end{figure}

    In addition to comparing the absolute computational cost of the wall models, it is also instructive to analyze how the wall models scale in the parallel environment. Figure \ref{fig:aposteriori_par_GQ} shows \textit{a posteriori} parallel performance analysis for GQ versus FV method through two metrics: speedup $S(p)$ and parallel efficiency $\eta(p)$ as a function of number of processors $p$ (or equivalently the number of control volumes per processor). The speedup $S(p)$ in figure \ref{fig:speedup_GQ_aposteriori} is defined as follows,

    \begin{equation}\label{eqn:speedup}
        S(p) = \frac{T_{p_{ref}}}{T_{p}},
    \end{equation}

    \noindent where $T_{p_{ref}}$ is the reference simulation time with $p_{ref}$ processors, $T_{p}$ is the simulation time with $p$ processors. Note that for each setup of wall model (no-slip, GQ or FV), $p_{ref}$ is defined with respect to that setup. In our tests we set $p_{ref}=32$, since a $p_{ref}$ below that overloaded the processors, causing the memory limit per processor to be exceeded. The ideal speedup follows the linear dashed line in figure \ref{fig:speedup_GQ_aposteriori}, that is given by $S(p)=p/p_{ref}$. From the figure, we observe sub-linear speedup for all the cases as the number of processors is increased. Interestingly, the speedup is sub-linear even for the no-slip case; whereas for the two wall models (GQWM and FVWM) the speedup is even lower than the no-slip case. To explain the trends in figure \ref{fig:aposteriori_par_GQ}, it is imperative to look into the issue of load imbalance in detail first.

\begin{figure}[t]
    	\centering
    	\begin{subfigure}[b]{0.46\textwidth}
    		\includegraphics[trim=80 220 100 220,clip,width=\textwidth]{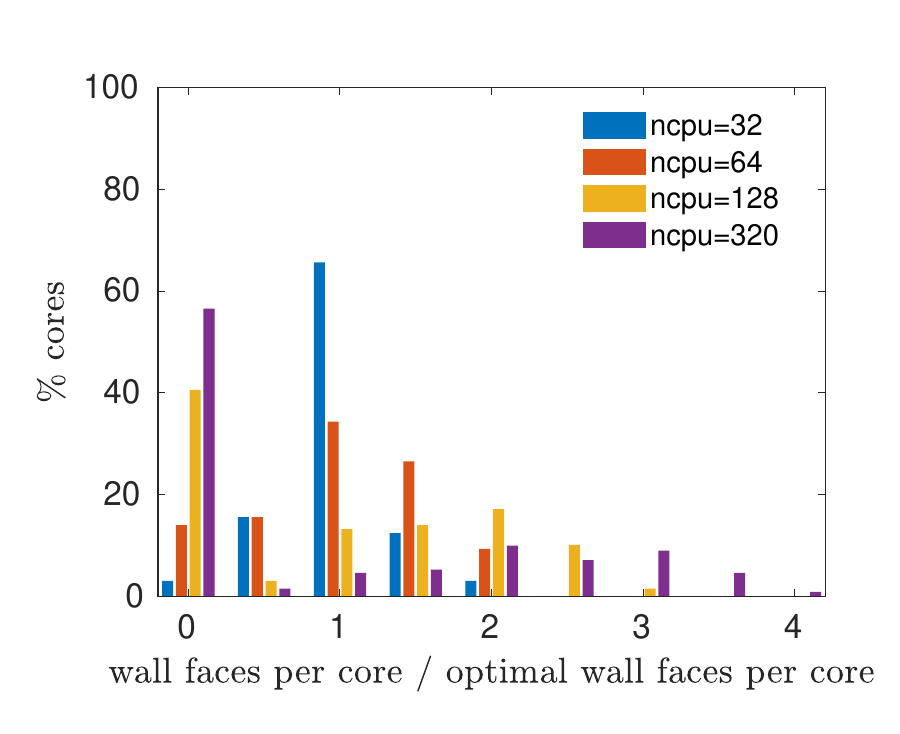}
    		\caption{}
    		\label{fig:load_imbalance_nfa}
    	\end{subfigure}
    	~
    	\begin{subfigure}[b]{0.46\textwidth}
    		\includegraphics[trim=80 220 100 220,clip,width=\textwidth]{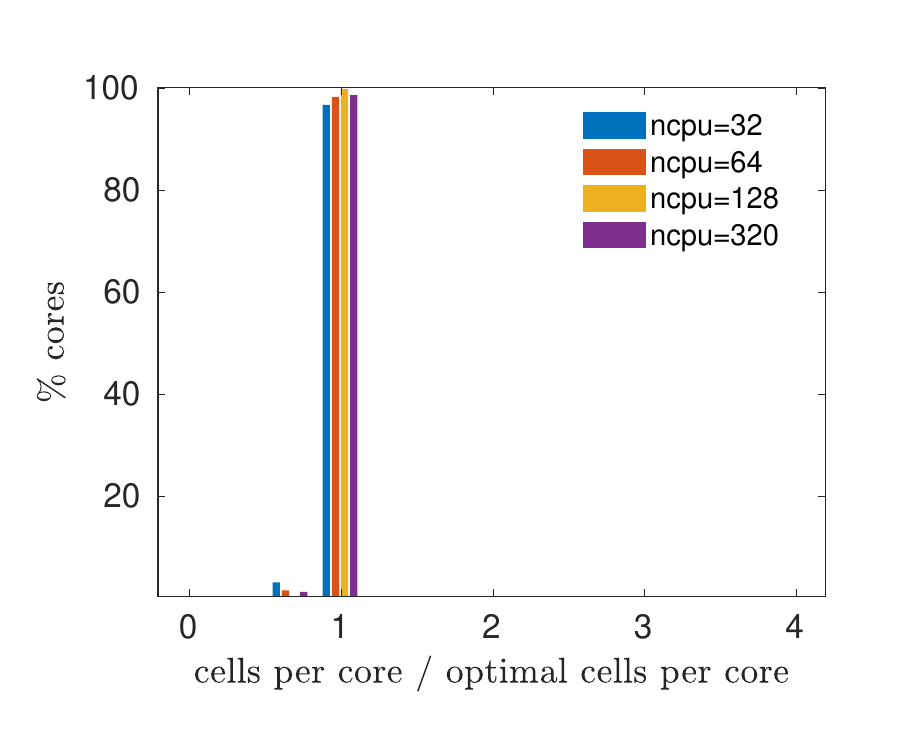}
    		\caption{}
    		\label{fig:load_imbalance_ncv}
    	\end{subfigure}    	
    	
    	\caption{\small Load imbalance among processors; histogram of resource distribution with respect to the load ((a) number of wall faces, (b) number of cells) taken by each processor. The horizontal axis is normalized by the optimal number of wall faces or cells per processor, which is given by the total number of wall faces or cells in the grid divided by the total number of processors used in the simulation.}
    	\label{fig:load_imbalance}
    \end{figure}

    When multiple processors are used in the simulation, the LES primal grid is partitioned among the processors in such a way that the computational load due to a particular aspect of the code (e.g. distribution of species, number of control volumes (ncv) or number of wall faces (nwfa) etc.) is aimed to be balanced among the processors. Here we make a distinction between the two types of load balances relevant to the current study, namely, the load balance in the LES solver (distribution of ncv) and the load balance in the wall-model solver (distribution of nwfa). In CharLES, the partitioning is done to optimize the load balance in the LES solver. However, this will likely result in a skewed distribution of wall faces among different processors, which results in load imbalance for the wall-model solver. This problem is elucidated in figure \ref{fig:load_imbalance} by plotting the number of wall faces in a processor against the percentage of processors containing those faces. For lower $p$, most of the processors have a loading close to the optimal loading of 1. As $p$ is increased, we observe an increasing percentage of processors containing no wall faces i.e. they remain idle during the wall-model solution step. Furthermore, a small percentage of processors are overly loaded with wall faces. To make this more concrete, for $p=320$, we observe that more than $50\%$ of the processors are idle in the wall-model solver whereas about $1\%$ of the processors are loaded with four times the ideal loading. It is clear from this figure, that by prioritizing only the load balance for LES solver, the load imbalance for the wall-model solver is exacerbated. Therefore, in WMLES, load imbalance due to wall-model solver can prove to be a much more serious issue than that due to LES solver.
    
    Referring back to figure \ref{fig:speedup_GQ_aposteriori}, the increasing departure of the no-slip line from the linear speedup with increase in $p$ is attributed to the increasing load imbalance in the LES solver. Note that although the load balance for LES solver is aimed to be optimized in CharLES, nevertheless it is far from ideal. For the two wall models, the additional drop in speedup compared to the no-slip simulation is attributed to the load imbalance in the wall-model solver. Note that the degradation in speedup is almost negligible for GQ method whereas FV method shows significant drop in speedup. Furthermore, the speedup degrades even further for FV method with more number of points $n$ used in the wall model. This is expected because with increase in $n$, the effect of load imbalance is now compounded with the higher computational cost of inverting large $n \times n$ matrices by the wall-model solver. More precisely, the difference in number of operations between the processors assigned with the most and the least number of wall faces is more for a higher $n$ than for a lower $n$.

    Figure \ref{fig:pareff_GQ_aposteriori} shows the parallel efficiency ($\eta$) for the three cases defined as,
    
    \begin{equation}\label{eqn:par_efficiency}
        \eta(p) = \frac{T_{p_{ref}} / (p/p_{ref})}{T_{p}},
    \end{equation}    

    \noindent where $\eta=1$ represents the ideal parallel scaling of a method. As expected from the speedup analysis, there is a parallel efficiency degradation for all three cases. However, $\eta$ for GQ method decays almost identically to the no-slip case, whereas FV method has a significant degradation of parallel efficiency, which aggravates further with increased number of points $n$.

    The above analysis indicates that FV method is more susceptible to adverse effects of load imbalance in the wall-model solver, whereas GQ method is remarkably agnostic to these effects. This is explained well by the observations made in \textit{a priori} tests regarding the operations count and the cost scaling versus $n$ for the two methods, that is, the GQ method circumvents the load imbalance because of the lower operations count associated with $n$ function evaluations compared to $n^2$ operation count for $n \times n$ matrix inversion in FV. A closer inspection also reveals that the convergence behavior with respect to iterations is quite rapid for GQ method, requiring only $2-3$ iterations versus $7-10$ iterations to convergence for the FV method. In summary, to alleviate the parallel-efficiency degradation, the FV method would require the additional implementation overhead of load balancing among the processors whereas GQ method circumvents this problem owing to its significantly lower operations count. Furthermore, as the GQWM cost is virtually insensitive to increase in number of quadrature points, the method would allow for an improvement in wall-stress accuracy without additional cost overhead.

    \subsection{\textit{A posteriori} performance characteristics of integral NEQWM}

    Figure \ref{fig:aposteriori_WMcost} shows that the computational cost of integral NEQWM LES is greater than that of the no-slip LES, but comparable to that of FVWM. In fact, for $p$ up to $384$, the additional cost of integral NEQWM is lower than that of FVWM with nominal number of grid points (i.e. $n=40$), as seen in figure \ref{fig:aposteriori_add_WMcost}. This is despite the fact that integral NEQWM is anticipated to represent the non-equilibrium effects in complex flows much better than the FV based ODEWM. Figure \ref{fig:aposteriori_par_GQ} shows that the parallel efficiency degradation of integral NEQWM with $p$ is quite drastic compared to FVWM. This is attributed to the communication overhead associated with {\tt MPI\_Allgatherv()} dominating the computational cost in the wall-model solver as $p$ is increased, because the wall-model solver must first collect wall quantities from all processors to fill the global wall-face array and then broadcast this array to all the cells in LES solver scattered across multiple processors, as described in \S \ref{sec:coupled_WMLES}, resulting in an \textit{all-to-all} communication. To alleviate the severe parallel efficiency degradation, a possible remedy could be to localize the broadcasting of wall quantities to the close vicinity of the wall, as the wall model takes gradients only from the cells adjacent to the wall. This will reduce the communication to relatively fewer partitions which lie close to the wall faces. However, implementing such a communication protocol could require additional implementation overhead, and therefore was not explored in the current study.

	\section{Conclusions}\label{sec:conclusion}
    
    Two wall models of integral form have been implemented and analyzed for cost efficiency and implementation overheads in an unstructured-grid finite-volume LES solver. Both the wall models avoid numerical integration in the wall-normal direction, thus significantly reducing the cost of the wall model compared with traditional approaches. The integral NEQWM solver presents unique implementation challenges with respect to data exchange with the LES solver and in terms of coordinate transformation of the wall quantities calculated within the wall model. The former challenge was resolved in the current study by using a one-way mapping of quantities from the wall faces to the cells in the LES solver. \textit{A priori} computational cost of the integral NEQWM is found to be comparable to that of the traditional ODE EQWM; however, the model suffers from severe parallel efficiency degradation due to communication overhead.
    
    For the Gauss-quadrature based ODEWM, a detailed computational cost and parallel efficiency analysis has been performed in both \textit{a priori} and \textit{a posteriori} settings. \textit{A priori} results highlight the favorable cost-scaling of the GQ method over the traditional FV method, with respect to the number of points used in the wall model, for a serial setting. The effects observed in \textit{a priori} tests have a direct bearing on the cost scaling in \textit{a posteriori} test, where the simulations are parallelized. As the number of processors are increased, we observe a perfect scaling for the GQWM (i.e. identical scaling to the no-slip simulations) indicating that GQ method does not contribute to the additional parallel-efficiency degradation due to the wall model. In contrast, FVWM encounters significant parallel-efficiency degradation, which can be attributed to the load-imbalance in the wall-model solver. In summary, the GQ method has two advantages over FV method: first, GQ method circumvents additional implementation overhead necessary to alleviate the parallel-efficiency degradation in FV owing to its significantly lower operations count; Second, as the GQWM cost is insensitive to the increase in number of quadrature points within a practical range, the method allows for an improvement of the predicted wall-stress accuracy without additional cost overhead.

    \section*{Acknowledgement}
    
    This investigation was funded by NASA (Grant 80NSSC18M0155) and the Office of Naval Research (ONR), Grant N00014- 16-S-BA10. Computing resources supporting this work were provided by the NASA High-End Computing (HEC) Program through the NASA Advanced Supercomputing (NAS) Division at Ames Research Center.
    
    \appendix
    \section{Near-wall asymptotic behavior of velocity for 3D integral NEQWM formulation}\label{appendix:A}
    
    To show the inconsistency in the original integral NEQWM formulation, we start with the linear term in the Taylor series expansion and use equation \ref{eqn:sublayer_inconsistent} to expand the term as follows,
    
    \begin{equation}\label{eqn:A1}
    u = \left.\frac{\partial u}{\partial y}\right|_{w} y = \frac{u_{\tau,x}}{\delta_{\nu}} y.
    \end{equation}  

    \noindent The viscous length scale $\delta_{\nu}$ is given by equation (C9) in \citep{yang15} as,
    
    \begin{equation}\label{eqn:A2}
    \delta_{\nu} = \nu \left[\frac{u_{\tau,x}^2 + u_{\tau,z}^2}{u_{\tau,x}^4 + u_{\tau,z}^4 }\right] = \frac{\mu}{\tau_{w}}\sqrt{\frac{\tau_{w,x}+\tau_{w,z}}{\rho}}.
    \end{equation}     

    \noindent Substituting \ref{eqn:A2} and definition of $u_{\tau,x} \equiv \left(\frac{\tau_{w,x}}{\rho}\right)^{1/2}$ in \ref{eqn:A1} we get the inconsistent asymptotic behavior,
    
    \begin{equation}\label{eqn:A3}
    u =\frac{\tau_{w}}{\mu}\sqrt{\frac{\tau_{w,x}}{\tau_{w,x}+\tau_{w,z}}} y = \frac{1}{\mu} \frac{\tau_{w,x}}{\cos{\theta}}\frac{1}{\sqrt{1+\tan{\theta}}} y= \left(\frac{\tau_{w,x}}{\mu}\right) \frac{1}{\sqrt{\cos^2\theta+\cos{\theta}\sin{\theta}}} y.
    \end{equation}      

    With the modified profile \ref{eqn:sublayer_modified}, the linear term in the Taylor series expansion and the viscous length scale $\delta_{\nu}$ are given by,
    
    \begin{equation}\label{eqn:A4}
    u = \left.\frac{\partial u}{\partial y}\right|_{w} y = \frac{u_{\tau,x}^2}{u_{\tau}\delta_{\nu}} y,
    \end{equation}     

    \begin{equation}\label{eqn:A5}
    \delta_{\nu} = \frac{\nu}{u_{\tau}},
    \end{equation}      

    \noindent Note that equation \ref{eqn:A5} is not a definition of $\delta_{\nu}$, rather it is determined from a consistent relation similar to equation (C9) in \citep{yang15}. Substituting \ref{eqn:A5} and definition of $u_{\tau,x} \equiv \left(\frac{\tau_{w,x}}{\rho}\right)^{1/2}$ in \ref{eqn:A4} we get,    

    \begin{equation}\label{eqn:A6}
    u = \left(\frac{\tau_{w,x}}{\mu}\right) y,
    \end{equation}      

    \noindent which is the desired consistent near-wall limiting behavior for velocity.
	
	\bibliography{mybibfile}
	
\end{document}